\documentclass[11pt, twoside]{article}
\usepackage{amsmath,amsthm}
\usepackage{amsfonts}
\usepackage{amssymb,latexsym}
\usepackage{enumerate}
\newtheorem{theorem}{Theorem}[section]
\newtheorem{remark}[theorem]{Remark}
\newtheorem{proposition}[theorem]{Proposition}
\newtheorem{lemma}[theorem]{Lemma}
\newtheorem{definition}[theorem]{Definition}

\newtheorem{corollary}[theorem]{Corollary}
\numberwithin{equation}{section}
\def\ind{1{\hskip -3 pt}\hbox{\textsc{I}}}
\def\n{\noindent}

\def\Ga{\Gamma}

\def\va{\varepsilon}

\def\o{\omega}
\def\O{\Omega}

\def\w{\wedge}

\def\om{\omega}

\def\ga{\gamma}
\def\wed{\wedge}
\def\Om{\Omega}

\newcommand{\SH}{{\rm SH_m}}

\newcommand{\MSH}{{\rm MSH}}

\numberwithin{equation}{section}

\begin{document}
	\setlength{\baselineskip}{18truept}
	\pagestyle{myheadings}

	\title {Degenerate complex Hessian equations with arbitrary measure in bounded domains}
	\author{
		 Nguyen Van Phu* and Le Mau Hai** \\
	* Faculty of Natural Sciences, Electric Power University,\\ Hanoi, Vietnam.\\
		**Department of Mathematics, Hanoi National University of Education ,\\ Hanoi,Vietnam;
		\\E-mail: phunv@epu.edu.vn    and  mauhai@hnue.edu.vn
	}

	\date{}
	\maketitle
	
	\renewcommand{\thefootnote}{}
	
	
	
	\renewcommand{\thefootnote}{\arabic{footnote}}
	\setcounter{footnote}{0}

\begin{abstract}
	\n
Let $\Omega$ be a bounded strictly $m$-pseudoconvex domain of $\mathbb{C}^n$. We solve degenerate complex Hessian equations of the form $(\omega + dd^c \varphi)^m\wedge\beta^{n-m} = \mu$ in the generalized Cegrell classes $\mathcal{K}_m(\Omega,\omega,\phi)$, where $\phi \in \mathcal{E}_m(\Omega)$ is a $m$-maximal function, $\omega$ is a smooth real $(1,1)$-form defined in a neighborhood of $\bar\Omega$ and $\mu$ is a positive Radon measure which is dominated by a Hessian measure of $m$-subharnomic functions in Cegrell class. 
\end{abstract}

	\section{Introduction}

 In 1982, E. Bedford and B. A. Taylor \cite{BT82} proved that the complex Monge-Amp\`ere operator $(dd^c.)^n$  is well-defined in the class of locally bounded plurisubharmonic functions over a bounded domain $\Omega$ in  $\mathbb{C}^n$ and whose image lies in the class of non-negative Borel measures. Given a positive Radon measure $\mu$ on $\Omega$,  in 1995,  S. Ko{\l}odziej  \cite{Ko95} established the existence of bounded solutions to the Monge-Amp\`ere equation $\mu=(dd^cu)^n$ under the assumption that  $\mu$ is dominated by  Monge-Amp\`ere measure of a bounded plurisubharmonic function. This result is widely known as Kołodziej’s subsolution theorem. Later on, U. Cegrell \cite{Ceg98, Ceg04} defined classes of functions  $\mathcal{F}(\Omega)$ and $\mathcal{E}(\Omega)$ which are not necessarily locally bounded, on which the complex Monge-Ampère operator is well-defined and continuous with respect to decreasing sequences. Based on this result, P. Ahag, U. Cegrell, R. Czyż, and P. H. Hiep \cite{ACCH09} proved that the Monge-Amp\`ere equation $\mu=(dd^cu)^n$ has solution in the class $\mathcal{E}(\Omega)$ if it has subsolution in this class. Subsequently,  L. M. Hai, P. H. Hiep, N. X. Hong and N. V Phu \cite[Theorem 3.1]{HHHP14} proved the subsolution theorem for the weighted complex Monge-Amp\`epe equation $\mu=-\chi(u)(dd^cu)^n$ in the class $\mathcal{E}_{\chi}(\Omega).$  
 
  Let $1\leq m\leq n$ be an integer. B{\l}ocki \cite{Bl05} and  Sadullaev, Abdullaev \cite{SA12} introduced and investigated the class of $m$-subharmonic functions which are a generalization of plurisubharmonic functions. They also proved that the complex $m$-Hessian operator $H_m(.) = (dd^c.)^m\wedge \beta^{n-m}$ is well-defined in the class of locally bounded $m$-subharmonic functions, where $\beta= dd^c\|z\|^2$ is the canonical K\"ahler form of $\mathbb{C}^n$. In \cite{Cu13}, N. N. Cuong proved the complex $m$-Hessian equation $\mu=H_m(u)$ has  bounded solutions in the class of $m$-subharmonic functions if it has a subsolution in this class. Building on Cegrell’s work, Lu Hoang Chinh \cite{Ch12,Ch15} introduced the classes $\mathcal{F}_m(\Omega)$ and $\mathcal{E}_m(\Omega)$ which are not necessarily locally bounded and  the complex $m$-Hessian operator  is well defined and continuous with respect to decreasing sequence. Subsequently, Nguyen Van Phu and coauthors, in a series of papers, studied subsolution theorems for functions in the Cegrell classes (see \cite{HP17,PD2023,PDtaiwan}).
  
  Let $\omega$ be a smooth real $(1,1)$-form defined in a neighborhood of $\bar\Omega$. Following the work of Bedford-Taylor \cite{BT76}, Ko\l odziej-Nguyen \cite{KN15Phong} defined the operator $(\omega + dd^c .)^n$ on locally bounded $\omega$-plurisubharmonic functions. They also proved that  the complex Monge-Ampère equation 
  \begin{equation}\label{Main eq}
  	(\omega + dd^c u)^n = \mu
  \end{equation}
  has a bounded  solution in the class of $\omega$-psh functions  if measure $\mu$ is  dominated by the Monge-Ampère measure of a bounded psh function (see \cite[Theorem 3.1]{KN23a}  and \cite[Theorem 3.2]{KN23b}). After that,  Salouf  \cite{Sal25} introduced the Cegrell-type classes $\mathcal{K}(\Omega,\omega, \phi),$ which are not necessarily bounded, and on which the operator $(\omega + dd^c .)^n$ is well defined and continuous along decreasing sequences.  Moreover,  when a maximal plurisubharmonic function $\phi$ is continuous on $\bar\Omega$ and the measure $\mu$ vanishes on pluripolar sets,  Salouf further  proved that the complex Monge-Ampère equation \eqref{Main eq} admits a solution in the class $\mathcal{K}(\Omega,\omega,\phi)$ provided that measure $\mu$ is dominated by the Monge-Ampère measure of a function in $\mathcal{K}(\Omega)$ \cite[Theorem 5.8]{Sal25}. 
  
    Recently, Ko\l odziej-Nguyen \cite{KN23c} and Gu-Nguyen \cite{GN18} defined the operator $(\omega + dd^c .)^m\w\beta^{n-m}$ on  locally bounded $(\omega,m)$-subharmonic functions. Moreover, Ko\l odziej-Nguyen \cite[Lemma 9.3]{KN23c} proved that the complex Hessian equation 
  \begin{equation}\label{Hessian eq}
  	(\omega + dd^c u)^m\w\beta^{n-m} = \mu
  \end{equation}
  has bounded solutions in the class of  $(\omega,m)$-subharmonic functions if measure $\mu$ is dominated by a Hessian measure of a bounded $m$-subharmonic function. In this paper, inspired by the work of  Salouf \cite{Sal25},  we defined the Cegrell-type class $\mathcal{E}_m(\Omega,\omega, \phi)$ which is not necessary bounded and on which the Hessian operator 	$(\omega + dd^c u)^m\w\beta^{n-m}$ is well-defined.  After that, we will prove that equation \eqref{Hessian eq} has unbounded solution in the class of  $(\omega,m)$-subharmonic functions if measure $\mu$ is dominated by a Hessian measure of a unbounded $m$-subharmonic function.\\
  
  The article is organized as follows. In Section 2, we recall some basic facts about Cegrell classes. In Section 3, we establish the Hessian operator on the class $\mathcal{E}_m(\Omega,\omega).$ In Section 4, we study  envelopes of $(\omega,m)$-subharmonic functions. The main result is proved in Section 5. In Section 6, we will prove a comparison principle.

\n 
{\bf Acknowledgments} 
The authors would like to thank Professor Nguyen Quang Dieu for valuable comments during the preparation of this work.\\
Some parts of this paper is written in our visit in Vietnam Institute for Advanced Study in Mathematics (VIASM) during the Spring of 2025. We also thank VIASM for its hospitality.

\section{Preliminaries} 
Throughout this paper, we will denote by $\Omega$ a bounded strictly $m$-pseudoconvex domain of $\mathbb{C}^n$
and denote by $\SH(\Omega)$ the set of $m$-subharmonic functions on $\Omega$.
Firstly, we recall  Cegrell classes introduced and investigated in  \cite{Ch12,Ch15}.

$$ \mathcal{E}^0_m(\Omega) = \{ u \in \SH(\Omega)\cap L^{\infty}(\Omega) : u = 0 \; \text{on} \; \partial{\Omega} \; \text{and} \; \int_\Omega (dd^c u)^m\w\beta^{n-m} < +\infty    \},
$$
$$ \mathcal{F}_m(\Omega) = \{ u \in \SH(\Omega) : \exists u_j \in \mathcal{E}^0_m(\Omega), \; u_j \searrow u \; \text{and} \; \sup_j \int_\Omega (dd^c u_j)^m\w\beta^{n-m}< +\infty    \}, $$ 
$$ \mathcal{E}_m(\Omega) = \{ u \in \SH(\Omega) : \forall D \Subset \Omega, \exists u_D \in \mathcal{F}(\Omega) \; \text{such that} \; u_D = u \; \text{on} \; D \}. $$
\n Relying on the work of Cegrell \cite{Ceg98,Ceg04}, Chinh \cite{Ch12, Ch15} proved that  the Hessian operator $H_m(.)=(dd^c .)^m\w\beta^{n-m}$ is well defined on $\mathcal{E}_m(\Omega)$ and continuous with respect to decreasing sequences. Moreover,  $\mathcal{E}_m(\Omega)$ is the largest set of negative $m$-sh functions for which these two above properties hold. Note that  a function in the class $\mathcal{E}_m(\Omega)$ is not necessary locally bounded and $\mathcal{E}_m(\Omega)$ contains all the $m$-sh functions that are  locally bounded and negative on $\Omega$. \\
Recall that a function $\phi \in \mathcal{E}_m(\Omega)$ is $m$-maximal if and only if it satisfies $(dd^c\phi)^m\w\beta^{n-m}=0.$ By $\MSH_m(\Omega),$ we denote the set of $m$-maximal functions in $\Omega$.\\
In \cite{T19}, Thien defined the class $\mathcal{N}_m(\Omega)$ as the set of functions $u \in \mathcal{E}_m(\Omega)$ for which the smallest $m$-maximal function 
majorant of $u$ is identical to zero.  In other words, 
$$ \mathcal{N}_m(\Omega) = \{ u \in \mathcal{E}_m(\Omega) : \forall \phi \in \mathcal{E}_m(\Omega)\cap \MSH_m(\Omega), \; u \leq \phi \leq 0 \Rightarrow \phi = 0 \}.   $$
We have 
$$ \mathcal{E}_m^0(\Omega) \subset \mathcal{F}_m(\Omega)  \subset \mathcal{N}_m(\Omega) \subset \mathcal{E}_m(\Omega). $$
If $\mathcal{K}_m(\Omega)\in\{\mathcal{E}_m(\Omega),\mathcal{N}_m(\Omega), \mathcal{F}_m(\Omega) \}$ and $\phi\in\mathcal{E}_m(\Omega),$ we say that a $m-$subharmonic $u$ defined on $\Omega$ belongs to $\mathcal{K}_m(\Omega,\phi)$ if there exists a function $\varphi\in\mathcal{K}_m(\Omega)$ such that $\phi\geq u\geq \varphi +\phi.$

Let $\omega$ be a smooth real $(1,1)$-form defined in a neighborhood of $\bar\Omega$. As in \cite{GN18} and \cite{KN23c}, we say that a function    $u:\Om \to [-\infty,+\infty)$ is called $(\omega,m)-\beta$-subharmonic if $u$ is upper semi-continuous, $u\in L^1_{\rm loc}(\Om,\beta^n)$ and for any collection $\ga_1,...,\ga_{m-1} \in \Ga_m(\beta)$,  we have
$$
(\omega+dd^c u )\wed \ga_1 \wed \cdots \wed \ga_{m-1} \wed \beta^{n-m}\geq 0
$$
in the sense of currents, where the positive cone $\Ga_m(\beta)$, associated to $\beta$, is defined as 
$$
\Ga_m(\beta):= \{\ga: \, \ga\ \ \text{is a real (1,1)}-\text{form}, \ga^k \wed \beta^{n-k} (z)>0  \quad\text{for } k=1,...,m, z\in \Om\}.
$$
For a $C^2$ function $u$, the above inequality is equivalent to inequalities

$$
(\o+dd^cu)^k\w\beta^{n-k}\geq 0\,\,\text{for}\,\,k=1,\ldots,m.	
$$
For convenience of presentation, in this paper we refer to an 
$(\omega,m)-\beta$-subharmonic  function simply as an 
$(\omega,m)$-subharmonic function.
The cone of $(\omega,m)$-subharmonic (resp. negative) functions on $\Om$ is denoted by $SH_{m} (\Om,\om)$ (resp. $SH_{ m}^{-} (\Om,\om)$).\\
Let $u\in SH_{m} (\Om,\om).$ We fix $\rho \in \SH(\Omega) \cap \mathcal{C}^\infty(\bar\Omega)$ such that $\rho = 0$ on $\partial \Omega$ and $\omega \leq dd^c \rho$. Then $u+\rho$ is a $m$-subharmonic function on $\O.$ We write $\tau=dd^c\rho-\o$ which is a smooth $(1,1)$-form, then $\o+dd^cu=dd^c(u+\rho)-\tau.$ Following the work of Bedford and Taylor \cite{BT82}, Ko\l odziej-Nguyen \cite[Section 9]{KN23c} and Gu-Nguyen \cite{GN18} prove that if $u\in SH_{m} (\Om,\om)\cap L^{\infty}_{loc}(\O)$ then the complex Hessian operator $$H_{m,\omega}(u)=(\omega+dd^cu)^m\w\beta^{n-m}:=\sum_{k=0}^m\binom{m}{k}(-1)^{m-k}[dd^c(u+\rho)^k]\w\tau^{m-k}\w\beta^{n-m}$$ is well-defined.\\

Note that, S. Ko{\l}odziej and N.-C. Nguyen \cite{KN23c} defined the Hessian operator $(\omega+dd^cu)^m\w\alpha^{n-m}$ for bounded $(\omega,m,\alpha)$-subharmonic function. So, it is natural and desirable to extend the Hessian operator $(\omega+dd^cu)^m\w\alpha^{n-m}$ for unbounded functions. However, we give a counterexample which shows that one cannot define $(dd^c u)^2 \wedge \alpha$ in a reasonably way  for arbitrary unbounded functions $u$ where $\omega=0, m=2$ and $n=3$. Indeed, consider the ball $\mathbb{B}=\mathbb{B}(0,\frac{1}{4})=\{(z_1,z_2,z_3)\in\mathbb{C}^3: |z_1|^2+|z_2|^2 + |z_3|^2<\frac{1}{4}\}$ in $\mathbb{C}^3,$  $\alpha=\frac{i}{2} dz_3\wedge d\bar{z}_3$ and consider the  function $u(z_1,z_2)= \sqrt{-\log|z_1|}(|z_2|^2-1)$ on $\mathbb{B}$. Then, by direct computation, we have
\begin{align*}&(dd^c u)^2\wedge\alpha = dd^cu\wedge dd^cu \wedge\alpha\\ 
	&=\bigg(2i\sum^3_{i,j=1}\frac{\partial^2u}{\partial z_i\partial\bar{z}_j}dz_i\wedge d\bar{z}_j\bigg)\wedge \bigg(2i\sum^3_{i,j=1}\frac{\partial^2u}{\partial z_i\partial\bar{z}_j}dz_i\wedge d\bar{z}_j\bigg)\wedge \frac{i}{2} dz_3\wedge d\bar{z}_3 \\
	&=32\bigg(\frac{\partial^2u}{\partial z_1\partial\bar{z}_1}.\frac{\partial^2u}{\partial z_2\partial\bar{z}_2}-\frac{\partial^2u}{\partial z_1\partial\bar{z}_2}.\frac{\partial^2u}{\partial z_2\partial\bar{z}_1}\bigg)dV_{2n}\\
	&=32\bigg( \frac{|z_2|^2-1}{16|z_1|^2\log|z_1|\sqrt{-\log|z_1|}}.\sqrt{-\log|z_1|}-\frac{-z_2}{4z_1\sqrt{-\log|z_1|}}.\frac{-\bar{z}_2}{4\bar{z}_1 \sqrt{-\log|z_1|}}\bigg)dV_{2n}\\
	&=\bigg( \frac{-2(|z_2|^2-1)}{|z_1|^2(-\log|z_1|)}- \frac{2|z_2|^2}{|z_1|^2(-\log|z_1|)}\bigg)dV_{2n}\\
	&=\frac{4}{(-\log|z_1|)|z_1|^2}(\frac{1}{2}- |z_2|^2)dV_{2n},
\end{align*} where $dV_{2n}=\frac{i}{2}dz_1\wedge d\bar{z}_1\wedge \frac{i}{2}dz_2\wedge d\bar{z}_2\wedge \frac{i}{2}dz_3\wedge d\bar{z}_3.$\\
Hence, we have $(dd^c u)^2\wedge\alpha\geq 0.$ Obviously, we have $dd^c u\wedge\alpha^2= 0.$ Thus,  we obtain $u$ is a $(0,2,\alpha)$-subharmonic function. However, $(dd^c u)^2\wedge \alpha$ is not a Radon measure on $\mathbb{B}$. Indeed, take $\varepsilon <\frac{1}{8}$ and consider $K=\mathbb{\overline{B}}(0,\varepsilon)$. Then we have $\int\limits_{K}(dd^c u)^2\wedge\alpha=+\infty$ and we get the desired conclusion.\\
In the  next section, we will define the Hessian operator $(\omega+dd^cu)^m\w\beta^{n-m}$ for a class of unbounded $(\omega,m)$-subharmonic functions.

\section{Definition of the Hessian operator $(\omega + dd^c.)^m\w\beta^{n-m}$ on $\mathcal{E}_m(\Omega,\om)$.}

\subsection{The current $(dd^c .)^k\w\beta^{n-m}$}
The aim of this Subsection is to show that the current $(dd^c u)^k\w\beta^{n-m}$ is well defined for any $u \in \mathcal{E}_m(\Omega)$ and with $1\leq k\leq m$. Note that, in the case when $u\in SH_m(\Omega)\cap L^{\infty}(\Omega)$, the current $(dd^c u)^k\w\beta^{n-m}$ is well defined as in \cite[Section 2]{SA12}. Firstly, we recall the Theorem 3.14 in \cite{Ch12}.
\begin{proposition}\label{bd2}
	Let $u^p\in\mathcal{E}_m(\O), p=1,\cdots,m$ and $(g^p_j)_j\subset \mathcal{E}^0_m(\Omega)$ such that $g^p_j\searrow u^p,\forall p=1,\cdots,m$. Then the sequence of measures 
	
	$$dd^c g^1_j\wed dd^c g^2_j\w \ldots\wed dd^c g^m_j\wed\beta^{n-m}$$
	\n converges weakly to a positive Radon measure which does not depend on the choice of the sequence $(g^p_j)_j$. We then define $dd^c u^1\wed dd^c u^2\w \ldots\wed dd^c u^m\wed\beta^{n-m}$
	is this weak limit.
\end{proposition}
\begin{remark}\label{rm4.4}
	Note that, by Lemma 3.13 in \cite{Ch12}, we have $$\mathcal{C}^{\infty}_0\subset \mathcal{E}^0_m(\Omega)\cap C(\Omega)-\mathcal{E}^0_m(\Omega)\cap C(\Omega).$$ Therefore,
	to prove the above Proposition, the author proved that for each $h\in\mathcal{E}^0_m(\Omega),$ we have 
	$$\lim\limits_{j\to\infty}\int_{\Omega}hdd^c g^1_j\wed dd^c g^2_j\w \ldots\wed dd^c g^m_j\wed\beta^{n-m}=\int_{\Omega}hdd^c u^1\wed dd^c u^2\w \ldots\wed dd^c u^m\wed\beta^{n-m}.$$ 
	\end{remark}
\n Next, we will prove the following Proposition which is a generalization of Proposition \ref{bd2}.
\begin{proposition}\label{pro4.5}
	Assume that $u^p, u^p_j\in\mathcal{E}_m(\O), p=1,\cdots,m$ and $u^p_j\searrow u^p,\forall p=1,\cdots,m.$ Then 	$dd^c u^1_j\wed dd^c u^2_j\w \ldots\wed dd^c u^m_j\wed\beta^{n-m}$ tends weakly to $dd^c u^1\wed dd^c u^2\w \ldots\wed dd^c u^m\wed\beta^{n-m}$ as $j\to\infty.$
\end{proposition}
\begin{proof}
	By Theorem 3.1 in \cite{Ch15}, there exists a sequence $[g^p_j]\subset\mathcal{E}_{m}^{0}(\Omega)\cap C(\overline{\Omega})$  which decreases to $u^p$ on $\Omega$ for all $p=1,\cdots,m.$  We put $v^p_j=\max(g^p_j,u^p_j)$ for all $p=1,\cdots,m.$ According to Theorem 3.9 in \cite{Ch15}, we infer that $v^p_j\in\mathcal{E}^0_m(\Omega)$ for all $p=1,\cdots,m.$ Obviously, we have $v^p_j\searrow u^p$ for all $p=1,\cdots,m$ as $j\to\infty.$ Thus, by Remark \ref{rm4.4}, we have
		\begin{equation}\label{eq4.1}\lim\limits_{j\to\infty}\int_{\Omega}hdd^c v^1_j\wed dd^c v^2_j\w \ldots\wed dd^c v^m_j\wed\beta^{n-m}=\int_{\Omega}hdd^c u^1\wed dd^c u^2\w \ldots\wed dd^c u^m\wed\beta^{n-m}.
			\end{equation}
		Since every function in $\mathcal{E}_m(\Omega)$ is locally the restriction of a function in $\mathcal{F}_m(\Omega)$ (see Remark 3.6 in \cite{Ch12}) and local property of weak convergence, we can assume that $u^p,u^p_j\in\mathcal{F}_m(\Omega)$ for all $p=1,\cdots,m.$ \\
		On the other hand, by integration by parts for functions in $\mathcal{F}_m(\Omega)$ (see Theorem 3.16 in \cite{Ch15}), we have
	\begin{equation}
		\begin{split}\label{eq4.3}
		&\int_{\Omega}hdd^c v^1_j\wed dd^c v^2_j\w \ldots\wed dd^c v^m_j\wed\beta^{n-m}\\
		&=	\int_{\Omega}v^1_jdd^c h\wed dd^c v^2_j\w \ldots\wed dd^c v^m_j\wed\beta^{n-m}\\
		&\geq \int_{\Omega}u^1_jdd^c h\wed dd^c v^2_j\w \ldots\wed dd^c v^m_j\wed\beta^{n-m}\\
		&= \int_{\Omega}hdd^c u^1_j\wed dd^c v^2_j\w \ldots\wed dd^c v^m_j\wed\beta^{n-m}\\
		&...\\
		&\geq \int_{\Omega}hdd^c u^1_j\wed dd^c u^2_j\w \ldots\wed dd^c u^m_j\wed\beta^{n-m}.
		\end{split}
		\end{equation}
	Using the same argument, we also have
	\begin{equation}\label{eq4.4} \int_{\Omega}hdd^c u^1_j\wed dd^c u^2_j\w \ldots\wed dd^c u^m_j\wed\beta^{n-m}\geq  \int_{\Omega}hdd^c u^1\wed dd^c u^2\w \ldots\wed dd^c u^m\wed\beta^{n-m}.
		\end{equation}
	Combining equality \eqref{eq4.1}, inequality \eqref{eq4.3} and inequality \eqref{eq4.4}, we infer that
		\begin{equation*}\label{eq4.5}\lim\limits_{j\to\infty}\int_{\Omega}hdd^c u^1_j\wed dd^c u^2_j\w \ldots\wed dd^c u^m_j\wed\beta^{n-m}=\int_{\Omega}hdd^c u^1\wed dd^c u^2\w \ldots\wed dd^c u^m\wed\beta^{n-m}.
			\end{equation*}
		Note that $$\mathcal{C}^{\infty}_0(\Omega)\subset \mathcal{E}^0_m(\Omega)\cap C(\Omega)-\mathcal{E}^0_m(\Omega)\cap C(\Omega).$$
		Thus, for $\chi\in \mathcal{C}^{\infty}_0(\Omega),$ we have
		$$\lim\limits_{j\to\infty}\int_{\Omega}\chi dd^c u^1_j\wed dd^c u^2_j\w \ldots\wed dd^c u^m_j\wed\beta^{n-m}=\int_{\Omega}\chi dd^c u^1\wed dd^c u^2\w \ldots\wed dd^c u^m\wed\beta^{n-m}.$$ The proof is complete.
	\end{proof}
\n Now let $ u \in \mathcal{E}_m(\Omega)$, and let $\alpha$ be a smooth $(m-k,m-k)$-form with ${\rm supp}\alpha\Subset\Om$. Take a relatively compact domain $D$ such that ${\rm supp}\alpha\Subset D\Subset \Om$. By \cite[Proposition 3.4]{Sal25}, we can write
$$ \alpha = \sum f_j T_j, $$
where $f_j$ are $\mathbb{C}$-values smooth functions on $\bar{\Om}$ and $T_j = dd^c u_1^j \wedge ... \wedge dd^c u_{m-k}^j$, $u_1^j,\ldots, u_{m-k}^j\in PSH(\Omega)\cap \mathcal{C}^{\infty}(\overline{\Omega})\subset \mathcal{E}_m(\Om).$  By Proposition \ref{bd2}, we note that $(dd^c u)^k\w\beta^{n-m}\w T_j=(dd^c u)^k\w dd^c u_1^j \wedge ... \wedge dd^c u_{m-k}^j\w\beta^{n-m} $ is a positive Radon measure. Set $T= (dd^c u)^k\w\beta^{n-m}$ and define $T$ acts on $\alpha$ by the equality
\begin{equation}\label{eq3.4}\left\langle T,\alpha\right\rangle= \sum\limits_{j}f_j(dd^c u)^k\w\beta^{n-m}\w T_j.\end{equation}

\n Then we have
$$|\left\langle T,\alpha\right\rangle|\leq\sum\limits_{j}|f_j|\int\limits_{\bar{D}}|(dd^c u)^k\w\ T_j\w\beta^{n-m}|.$$
\noindent On the other hand, since $(dd^c u)^k\w\beta^{n-m} \wedge T_j $ is a positive Radon measure then $\int\limits_{\bar{D}}|(dd^c u)^k\w\ T_j\w\beta^{n-m}|$ is finite and, hence, the right-hand side of the above inequality is bounded by a constant $C,$ where $C$ is a constant dependent on $D,\Omega$. Thus by definition of currents in \cite{Kl91}, $T=(dd^c u)^k\w\beta^{n-m}$ defines a $(n+k-m, n+k-m)$-current. Such as,  we can define  $(dd^c u)^k\w\beta^{n-m}$ acts on a smooth $(m-k,m-k)$-form $\alpha$ by 
$$\left\langle (dd^c u)^k\w\beta^{n-m}, \alpha\right\rangle= \sum f_j (dd^c u)^k\w\beta^{n-m} \wedge T_j, $$

\noindent if $\alpha=  \sum\limits_{j\in J} f_j T_j $. \\
To make the above definition reasonable, we need to show that this definition does not depend on the representation of $\alpha.$ Indeed, suppose now that 
$$ \alpha = \sum f_j T_j = \sum g_l S_l, $$
for $f_j, g_l \in \mathcal{C}^\infty(\bar\Omega)$ and $T_j, S_l$ are as in equality \eqref{eq3.4}. We need to prove that 
$$ \sum f_j (dd^c u)^k\w \beta^{n-m}\wedge T_j = \sum g_l (dd^c u)^k\w\beta^{n-m} \wedge S_l. $$
 To get this, we put $u_s=\max(u,-s)$. Since $u_s\in SH_m(\Omega)\cap L^{\infty}(\Omega),$ according to \cite[Section 2]{SA12} we have the current $(dd^cu_s)\wedge\beta^{n-m}$ is well-defined. Therefore, we have
\begin{equation}\label{e4.5}(dd^cu_s)^k\w\beta^{n-m}\w\alpha=\sum f_j (dd^c u_s)^k\w\beta^{n-m} \wedge T_j = \sum g_l (dd^c u_s)^k\w\beta^{n-m} \wedge S_l, \quad \forall s.
	\end{equation}
Note that, we also have $u_s\in \mathcal{E}_m(\Omega)$  and  $u_s \searrow u$ on $\Omega$ as $s\to\infty$. According to Proposition \ref{pro4.5}, let $s\to\infty$,  we deduce that
$$(dd^c u_s)^k\w\beta^{n-m} \wedge T_j= (dd^c u_s)^k\w dd^c u_1^j \wedge ... \wedge dd^c u_{m-k}^j\w\beta^{n-m}$$ converges weakly to 
$$(dd^c u)^k\w\beta^{n-m} \wedge T_j= (dd^c u)^k\w dd^c u_1^j \wedge ... \wedge dd^c u_{m-k}^j\w\beta^{n-m}.$$
Since $f_j\in\mathcal{C}^{\infty}(\overline{\Omega})$, we infer that
$f_j(dd^c u_s)^k\w\beta^{n-m} \wedge T_j$ converges weakly to $f_j(dd^c u)^k\w\beta^{n-m} \wedge T_j$ as $s\to\infty.$\\
 Similarly, we also have $g_l (dd^c u_s)^k\w\beta^{n-m} \wedge S_l$ converges weakly to $g_l (dd^c u)^k\w\beta^{n-m} \wedge S_l$ as $s\to\infty.$\\
It follows from equation \eqref{e4.5} that
$$ \sum f_j (dd^c u)^k\w\beta^{n-m} \wedge T_j = \sum g_l (dd^c u)^k \w\beta^{n-m}\wedge S_l. $$ This proved that our definition is meaningful.\\
From the above result, we give the following theorem.
\begin{theorem}\label{thm 4.2}
	Let $ u \in \mathcal{E}_m(\Omega)$. For all $k \in \{1,..,m \}$, the current $(dd^c u)^k\w\beta^{n-m}$ is well defined. Furthermore, if  $ (u_j)_j $  is a decreasing sequence in $ \mathcal{E}_m(\Omega) $ that converges to $ u$, then the sequence $ \left( ( dd^c u_j)^k\w\beta^{n-m} \right)_j $ converges weakly  to $ ( dd^c u)^k\w\beta^{n-m}$.  
\end{theorem}

\subsection{The classes $\mathcal{E}_m(\Omega,\omega)$}
First, we have the following definition.
\begin{definition}
	By $\mathcal{P}_{m,\om}(\Om)$ we denote the set of $\rho\in C^2(\bar{\Omega})\cap \SH(\Om)$ such that $\rho=0$ on $\partial\Om$ and $\om\leq dd^c\rho$.
\end{definition}
\n Now, we give the definition of the class $\mathcal{E}_m(\Omega,\om).$
\begin{definition}
	We say that $u\in\mathcal{E}_m(\Om,\om)$ (resp. $\mathcal{N}_m(\Om,\om),\mathcal{F}_m(\Om,\om)$) if $u\in \SH(\Om,\om)$ and $u+\rho\in \mathcal{E}_m(\O)$ (resp. $\mathcal{N}_m(\O), \mathcal{F}_m(\O)$).
\end{definition}

\n From now to the end of the paper, we denote by $\phi$ a maximal $m$-sh function in $\mathcal{E}_m(\Omega)$. \\
Fix $\rho \in \mathcal{P}_{m,\omega}(\Omega)$. Let $\mathcal{K}_m(\Omega)\in\{\mathcal{E}_m(\Omega),\mathcal{N}_m(\Omega), \mathcal{F}_m(\Omega) \}.$ Similar to \cite[Subsection 4.2]{Sal25}
we define the set
$\mathcal{K}_m(\Omega,\omega, \phi)$ as follow
\begin{definition}\label{de3.7}
$$ u \in  \mathcal{K}_m(\Omega,\omega, \phi) \Leftrightarrow u \in \SH(\Omega,\omega) \; \text{and} \;  u + \rho \in \mathcal{K}_m(\Omega,\phi). $$
\end{definition}
\begin{remark}
{\rm	Note that the definition \ref{de3.7} does not depend on $\rho$.
	Indeed,
	let $\rho' \in \mathcal{P}_{m,\omega}(\Omega) $. If $u+\rho \in \mathcal{K}_m(\Omega,\phi)$, then we have 
	$$\phi + \tilde{u} \leq u+\rho  \leq \phi$$
	for some $\tilde{u}\in \mathcal{K}_m(\Omega)$.
	Since $\phi$ is maximal and by $\lim\sup\limits_{z\to x} (u+\rho')(z)\leq\phi(x)$ for all $x\in\partial{\Om}$ then  we get $u+\rho' \leq \phi$ on $\Om$. On the other hand 
	$$ u+\rho' \geq u+\rho+\rho' \geq \phi +\tilde{u}+\rho',$$
	where $\tilde{u}+\rho'\in \mathcal{K}_m(\Omega)$ and the desired conclusion follows.}
\end{remark}
\n We have the following observation:
\begin{proposition}
	
(i)	$\mathcal{K}_m(\Omega,\omega,\phi) \subset \mathcal{E}_m(\Omega,\omega). $\\
		(ii) $\mathcal{E}_m(\Omega,0) = \mathcal{E}_m(\Omega)$.
		
\end{proposition}
\begin{proof}
(i)	Let $\rho \in \mathcal{P}_{m,\omega}(\Omega)$. 
	If $u \in \mathcal{K}_m(\Omega,\omega,\phi)$, then $u+\rho \in \mathcal{K}_m(\Omega,\phi)$.  It follows that
	$$  0 \geq \phi \geq u+\rho \geq \phi + \tilde{u}, $$
	for some $\tilde{u} \in \mathcal{K}_m(\Omega)$. Hence
	$u \in \mathcal{E}_m(\Omega,\omega)$.\\
	(ii) We take $\rho = 0$ in the definition of $\mathcal{E}_m(\Omega,0)$ to get the second assertion.
\end{proof}
\n Next, we will show  that the set $\mathcal{E}_m(\Omega,\omega)$ is always a non-empty subset of $\SH(\Omega,\omega)$.
\begin{proposition}\label{pro3.10}
	$$ \SH^-(\Omega, \omega) \cap L^\infty_{loc}(\Omega) \subset \mathcal{E}_m(\Omega,\omega). 
	$$
\end{proposition}
\begin{proof}
	Let $ \rho \in \mathcal{P}_{m,\omega}(\Omega) $. If $ u \in \SH^-(\Omega,\omega)\cap L^\infty_{loc}(\Omega), $ then $ u + \rho \in \SH^-(\Omega) \cap L^\infty_{loc}(\Omega) \subset \mathcal{E}_m(\Omega) $  and therefore $ u \in \mathcal{E}_m(\Omega,\omega)$.
\end{proof}

\subsection{The Hessian operator $( \omega + dd^c.)^m\w\beta^{n-m}$ on $\mathcal{E}_m(\Omega,\omega)$} \phantom{} \\
Let $u \in \mathcal{E}_m(\Omega, \omega)$, and let $ \rho \in \mathcal{P}_{m,\omega}(\Omega)$. Write
\begin{align*}& ( \omega + dd^c u)^m\w\beta^{n-m} = \Big( \omega - dd^c \rho + dd^c ( u+ \rho) \Big)^m\w\beta^{n-m}\\
	& = \sum_{k=0}^m \binom{m}{k} (-1)^{m-k} ( dd^c (u+\rho))^k \wedge \tau^{m-k}\w\beta^{n-m}, 
	\end{align*}
where $ \tau = dd^c \rho - \omega $ is a smooth semi-positive real $(1,1)$-form.  Because $ u +\rho\in\mathcal{E}_m(\Om)$, then by the proof before  Theorem \ref{thm 4.2},  
$( dd^c (u+\rho) )^k \wedge \tau^{m-k} \w\beta^{n-m}$  defines a positive Radon measure. By linearity, we can define 
the operator $ (\omega + dd^c u)^m\w\beta^{n-m}$. We have to check that $(\omega+dd^c u)^m\w\beta^{n-m}$ is independent of $ \rho$. Let $ \rho' \in \mathcal{P}_{m,\omega}(\Omega)$. It remains then to prove that 
\begin{align*}
	& \sum_{k=0}^m \binom{m}{k} (-1)^{m-k} ( dd^c (u+\rho))^k \wedge (dd^c \rho - \omega)^{m-k}\w\beta^{n-m} \\ =  &\sum_{k=0}^m \binom{m}{k} (-1)^{m-k} ( dd^c (u + \rho'))^k \wedge (dd^c \rho' - \omega)^{m-k}\w\beta^{n-m}.
\end{align*}
Indeed, we put $u_j=\max(u,-j)$. We have $u_j+\rho\geq u+\rho\in\mathcal{E}_{m}(\Omega).$ Thus, we imply that $u_j\in\mathcal{E}_m(\Omega, \omega)\cap L^{\infty}(\Omega)$.   Now,
we have \begin{align*}
	&( \omega + dd^c u_j)^m\w\beta^{n-m}\\
	&= \sum_{k=0}^m \binom{m}{k} (-1)^{m-k} ( dd^c (u_j+\rho))^k \wedge (dd^c \rho - \omega)^{m-k}\w\beta^{n-m} \\ =  &\sum_{k=0}^m \binom{m}{k} (-1)^{m-k} ( dd^c (u_j + \rho'))^k \wedge (dd^c \rho' - \omega)^{m-k}\w\beta^{n-m}.
\end{align*}
 Letting $j \rightarrow +\infty$, the result follows from  Theorem \ref{thm 4.2}. 
\begin{definition}
	Let $ u \in \mathcal{E}_m(\Omega, \omega)$. Then the operator $ (\omega + dd^c u)^m \w\beta^{n-m}$ is given by the  formula
	\begin{equation}\label{eq4.6} ( \omega + dd^c u)^m\w\beta^{n-m} = \sum_{k=0}^m \binom{m}{k} (-1)^{m-k} ( dd^c (u+\rho))^k \wedge (dd^c \rho - \omega)^{m-k}\w\beta^{n-m},\end{equation}
	where $ \rho \in \mathcal{P}_{m,\omega}(\Omega)$.
\end{definition}
As a consequence of Theorem \ref{thm 4.2}, we obtain the following result.
\begin{theorem}\label{mono}
	The Monge-Ampère operator $ (\omega + dd^c .)^m\w\beta^{n-m} $ is well  defined on the class $\mathcal{E}_m(\Omega,\omega)$. Furthermore, if $ (u_j)_j \subset {\mathcal{E}}_m(\Omega,\omega)$ is a decreasing sequence that converges to $u \in \mathcal{E}_m(\Omega,\omega)$, then the sequence of Radon measures $\left( ( \omega + dd^c u_j )^m \w\beta^{n-m}\right)_j$ converges weakly to $ ( \omega + dd^c u)^m\w\beta^{n-m}$.
\end{theorem}

\begin{remark}
{\rm Let $u\in \mathcal{E}_m(\Omega,\omega).$ Thanks to \cite[Theorem 3.18]{GN18}, there exists $ (u_j)_j \subset {\mathcal{E}}_m(\Omega,\omega)\cap C^{\infty}(\O)$ which is a decreasing sequence and it converges to $u \in \mathcal{E}_m(\Omega,\omega).$ Obviously, we have $( \omega + dd^c u_j )^m \w\beta^{n-m}$ is a positive Radon measure. Thus, by Theorem \ref{mono}, we obtain that $ ( \omega + dd^c u)^m\w\beta^{n-m}$ is a positive Radon measure.}
	\end{remark}

	\begin{definition}\label{dn 5.1}
		We denote by $ \mathcal{K}^a_m(\Omega,\omega)$ the set of $ u \in \mathcal{K}_m(\Omega,\omega)$ such that $(\omega + dd^c u)^m\w\beta^{n-m}$ puts no mass on $m$-polar set.
	\end{definition}
\begin{definition}\label{dn 5.2}
		 By $ \mathcal{K}^a_m(\Omega,\omega,\phi)$ we denote the set of $ u \in \mathcal{K}_m(\Omega,\omega,\phi)$ such that there exists a function $\tilde{u}\in \mathcal{K}^a_m(\Omega)$ with $\tilde{u}+\phi\leq u+\rho\leq \phi,$ where $\tilde{u}\in \mathcal{K}^a_m(\Omega)$ means that $H_m(\tilde{u})$ vanishes on all $m$-polar sets.
	\end{definition}

	\begin{proposition}\label{prop 5.2}
		If  $ u \in \mathcal{F}^a_m(\Omega,\omega)\cup [SH_{m}(\Omega,\omega)\cap L^{\infty}(\Omega)]$, then  $u+\rho\in\mathcal{E}^a_m(\O)$ for all $\rho\in\mathcal{P}_{m,\o}.$ 
	\end{proposition}
	\begin{proof} Obviously, the above Proposition is true in the case when $u\in SH_{m}(\Omega,\omega)\cap L^{\infty}(\Omega).$ It remains to prove in the case when $u \in \mathcal{F}^a_m(\Omega,\omega).$ Indeed, let $A$ be a $m$-polar subset of $ \Omega$. We have 
		\begin{align*}
			 &0=\int_A (\omega + dd^c u)^m\w\beta^{n-m} 
			= \int_A (dd^c (u +\rho))^m\w\beta^{n-m}\\
			& + \sum_{k=1}^m (-1)^k \binom{m}{k}\int_A (dd^c \rho - \omega)^k \wedge (dd^c(u+\rho))^{m-k}\w\beta^{n-m}.
		\end{align*}
		Let $\sigma \in \mathcal{P}_{m,-\omega}(\Omega)$. \cite[Lemma 5.6]{HP17} gives for all $ 1 \leq k \leq m$:
		\begin{align*}
			&\int_A (dd^c\rho- \omega)^k \wedge (dd^c(u+\rho))^{m-k}\w\beta^{n-m} \\
			&\leq  \int_A (dd^c (\rho+\sigma)) )^k \wedge (dd^c(u+\rho))^{m-k}\w\beta^{n-m} \\
			&\leq \left( \int_A (dd^c (\rho+\sigma))^m\w\beta^{n-m} \right)^{k/m} \left( \int_A (dd^c(u+\rho))^m\w\beta^{n-m} \right)^{(m-k)/m} \\
			&= 0,
		\end{align*}
	where the last equality follows from the fact that $\rho+ \sigma$ is in $\SH(\Omega)\cap C^2(\bar{\Om})$ and $\int_A (dd^c(u+\rho))^m\w\beta^{n-m}\leq \int_{\Omega} (dd^c(u+\rho))^m\w\beta^{n-m}<+\infty$ because $u+\rho\in\mathcal{F}_m(\Omega) $ (see Theorem 4.9 in \cite{T19}).
		The proof is complete.
	\end{proof}
\begin{remark}\label{rm5.3}
{\rm (i) Because Lemma 5.6 in \cite{HP17} also holds for the class $\mathcal{E}_m(\Om)$ then repeating the proof of Proposition \ref{prop 5.2}, we see that if $u+\rho\in\mathcal{E}^a_m(\O)$ for some $\rho\in\mathcal{P}_{m,\o},$ then $(\omega + dd^c u)^m\w\beta^{n-m} $ vanishes on all $m$-polar sets.\\
(ii) 	For every $u \in \mathcal{E}_m(\Omega,\omega)$, we have 
$$ \ind_{\{u=-\infty\}}(\omega + dd^c u)^m\w\beta^{n-m} = \ind_{\{u=-\infty\}}\big(dd^c (u+\rho)\big)^m\w\beta^{n-m} . $$}
\end{remark}
 Using  \cite[Proposition 5.2]{HP17} and Remark \ref{rm5.3} (ii), we have the following proposition. 

\begin{proposition}\label{prop 2.5}
	If $u$, $v \in \mathcal{E}_m(\Omega,\omega)$ are such that $u \leq v$, then
	$$ \ind_{\{u=-\infty\}}(\omega + dd^c u)^m\w\beta^{n-m} \geq \ind_{\{v=-\infty\}}(\omega + dd^c v)^m\w\beta^{n-m}.$$
\end{proposition}
\n Next, we will prove maximum principle in the class $\mathcal{E}_m(\Omega,\omega).$
	\begin{theorem}\label{thm 5.3}
		If $u, v \in \mathcal{E}_m(\Omega,\omega)$, then 
		$$ \ind_{\{u>v\}} ( \omega + dd^c u)^m\w\beta^{n-m} = \ind_{\{u>v\}} ( \omega + dd^c \max(u,v))^m\w\beta^{n-m}. $$
	\end{theorem}
	\begin{proof}
		Fix $\rho \in \mathcal{P}_{m,\omega}(\Omega)$.  Write 
		$$ ( \omega + dd^c u)^m\w\beta^{n-m}  = \sum_{k=0}^m \binom{m}{k} (-1)^{m-k} ( dd^c (u+\rho))^k \wedge (dd^c \rho - \omega)^{m-k}\w\beta^{n-m}. $$
		\n On the other hand, by the hypothesis $u+\rho, v+\rho\in \mathcal{E}_m(\Om)$ then so is $\max(u+\rho, v+\rho)\in \mathcal{E}_m(\Om)$. Note that $\max(u+\rho, v+\rho)=\max(u,v)+\rho$ on $\Om$. Therefore, we  have $\rho+\max(u,v) \in\mathcal{E}_m(\Om)$. It means that $\max(u,v) \in\mathcal{E}_m(\Om,\o).$ Thus, the Hessian operator $(\om + dd^c\max(u,v))^m\w\beta^{n-m}$ is well defined.\\
		 Write 
		\begin{align*}& (\om + dd^c\max(u,v))^m\w\beta^{n-m} \\
			& = \sum_{k=0}^m \binom{m}{k} (-1)^{m-k} ( dd^c (\max(u,v)+\rho))^k \wedge (dd^c \rho - \omega)^{m-k}\w\beta^{n-m}. \\
			&= \sum_{k=0}^m \binom{m}{k} (-1)^{m-k} ( dd^c (\max(u+\rho,v+\rho)))^k \wedge (dd^c \rho - \omega)^{m-k}\w\beta^{n-m}
			\end{align*}
		
		\n Since $\omega$ is a smooth real $(1,1)$-form, we have $\tau=dd^c\rho-\omega$ is also a smooth real $(1,1)$-form. This implies that $\tau^{m-k}$ is a smooth real $(m-k,m-k)$-form. Repeating the arguments as in Proposition 4.1 in \cite{Sal25}, we can write
		$$ \tau^{m-k}= \sum\limits_{i\in I} f_i T_i, $$
		where $I$ is a finite set, $(f_i)_{i\in I}$ are smooth functions with real values  and $T_i=dd^cu^i_1\w\cdots\w dd^cu^i_{m-k}$  where, for every $i\in I$ and every $l=1,\cdots,m-k, u_l^i$ is a smooth negative plurisubharmonic function defined in a neighborhood of $\overline{\Omega}.$\\
		 By the linearity, it suffices to prove that
		\begin{equation}\label{eq5.1} \ind_{\{ u > v\}} (  dd^c \max(u+\rho,v+\rho))^k \wedge T_i\w\beta^{n-m} = \ind_{\{ u > v\}} ( dd^c (u+\rho))^k \wedge T_i\w\beta^{n-m}
			\end{equation}
		
		\n with $u+\rho, v+\rho\in \mathcal{E}_m(\Om)$, $T_i= dd^c u^i_1\w\cdots\w dd^c u^i_{m-k}$, where $u^1_1,\ldots, u^i_{m-k}\in PSH(\Omega)\cap \mathcal{C}^{\infty}(\overline{\Omega})\subset SH_m(\Omega)\cap L^{\infty}(\Omega)\subset \mathcal{E}_m(\Om).$ \\
		Because it is independent interest, we will prove it in the following Proposition.
		\begin{proposition} \label{pro5.6}
			Let  $u,u_1,\ldots,u_{m-k}\in \mathcal{E}_m(\Omega)$, $v\in\text{SH}^-_m(\Omega)$.  Then we have
			\begin{align*}&\ind_{\{u>v\}}[dd^c\max(u,v)]^k\wedge dd^cu_1\wedge \cdots\wedge dd^cu_{m-k}\wedge\beta^{n-m}\\
				&=\ind_{\{u>v\}}[dd^cu]^k\wedge dd^cu_1\wedge \cdots\wedge dd^cu_{m-k}\wedge\beta^{n-m}.
			\end{align*}
		\end{proposition}
	 \n We will use Theorem 3.6 in \cite{HP17} to prove Proposition \ref{pro5.6}. For the reader's convenience, we  recall Theorem 3.6 in  \cite{HP17}.
	\begin{lemma}\label{lm5.6}
		Let  $u,v_1,\ldots,v_{m-1}\in \mathcal{E}_m(\Omega)$, $v\in\text{SH}^-_m(\Omega)$.  Then
		\begin{align*}&\ind_{\{u>v\}}dd^c\max(u,v)\wedge dd^cv_1\wedge \cdots\wedge dd^cv_{m-1}\wedge\beta^{n-m}\\
			&=\ind_{\{u>v\}}dd^cu\wedge dd^cv_1\wedge \cdots\wedge dd^cv_{m-1}\wedge\beta^{n-m}.
		\end{align*}
	\end{lemma}
	\begin{proof}[The proof of Proposition \ref{pro5.6}]
		 Applying Lemma \ref{lm5.6} with $v_1=\max(u,v); $ we have
			\begin{equation}\begin{split}\label{eq5.2} &\ind_{\{ u > v\}}   [dd^c \max(u,v)]^2 \wedge dd^cv_2\wedge dd^cv_3\w\cdots\w dd^cv_{m-1}\w\beta^{n-m}\\
				& = \ind_{\{ u > v\}}  dd^c u\w dd^c \max(u,v) \wedge dd^cv_2\wedge dd^cv_3\w\cdots\w dd^cv_{m-1}\w\beta^{n-m}\\
				&= \ind_{\{ u > v\}}  dd^c \max(u,v)\wedge dd^cu \wedge dd^cv_2\wedge dd^cv_3\w\cdots\w dd^cv_{m-1}\w\beta^{n-m}\\
				&=\ind_{\{ u > v\}}  dd^cu\wedge dd^cu \wedge dd^cv_2\wedge dd^cv_3\w\cdots\w dd^cv_{m-1}\w\beta^{n-m}\\
				&=\ind_{\{ u > v\}}[dd^cu]^2\wedge dd^cv_2\wedge dd^cv_3\w\cdots\w dd^cv_{m-1}\w\beta^{n-m},
			\end{split}
				\end{equation}
			where the third equality is due to Lemma \ref{lm5.6} with $v_1=u.$\\
			Applying equality \eqref{eq5.2} with with $v_2=\max(u,v); $ we have
			\begin{equation}\begin{split}\label{eq5.3} &\ind_{\{ u > v\}}   [dd^c \max(u,v)]^3 \wedge dd^cv_3\wedge dd^cv_4\w\cdots\w dd^cv_{m-1}\w\beta^{n-m}\\
					& = \ind_{\{ u > v\}}  [dd^cu]^2\w dd^c \max(u,v) \wedge dd^cv_3\wedge dd^cv_4\w\cdots\w dd^cv_{m-1}\w\beta^{n-m}\\
					&= \ind_{\{ u > v\}}  dd^c \max(u,v)\wedge  [dd^cu]^2 \wedge dd^cv_3\wedge dd^cv_4\w\cdots\w dd^cv_{m-1}\w\beta^{n-m}\\
					&=\ind_{\{ u > v\}}  dd^cu\wedge [dd^cu]^2 \wedge dd^cv_3\wedge dd^cv_4\w\cdots\w dd^cv_{m-1}\w\beta^{n-m}\\
					&=\ind_{\{ u > v\}}[dd^cu]^3\wedge dd^cv_3\wedge dd^cv_4\w\cdots\w dd^cv_{m-1}\w\beta^{n-m},
				\end{split}
			\end{equation}
		where the third equality is due to Lemma \ref{lm5.6} with $v_1=v_2=u.$\\
\n By repeating this process, the desired conclusion follows. The proof of Proposition \ref{pro5.6} is complete.
	\end{proof}
Now, we continue the proof of Theorem \ref{thm 5.3}. Applying Proposition \ref{pro5.6} with with $u+\rho$ plays the role of $u;$ $v+\rho$ plays the role of $v;$ $T_i$ plays the role of $dd^cu_{1}\wedge\cdots\wedge dd^c u_{m-k}$  we obtain equality \eqref{eq5.1}. The proof of Theorem \ref{thm 5.3} is complete.
\end{proof}
\n The following Corollary is an extension of J. P. Demailly's result \cite[Proposition 6.11]{Dem91}, in which he proved this result when $m = n$ and $\omega = 0$.
	\begin{corollary}\label{cor 5.5}
		Let $u, v \in \mathcal{E}_m(\Omega,\omega)$. If  $(\omega + dd^c v)^m\w\beta^{n-m}$ 
		vanishes on $m$-polar sets, then 
		$$ ( \omega + dd^c \max(u,v))^m\w\beta^{n-m} \geq \textit{1}_{\{u> v \}} ( \omega + dd^c u)^m\w\beta^{n-m} + \textit{1}_{\{ u \leq v \}} ( \omega + dd^c v)^m\w\beta^{n-m}. $$
		In particular, if moreover $u \geq v$ then 
		$$ \textit{1}_{\{u=v\}} ( \omega + dd^c u)^m\w\beta^{n-m} \geq  \textit{1}_{\{u=v\}} ( \omega + dd^c v)^m\w\beta^{n-m}. $$
	\end{corollary}
	\begin{proof}
\n First, we note that if $u, v\in\mathcal{E}_m(\Om,\om)$ then $\max\{u,v\}\in\mathcal{E}_m(\Om,\om)$. Futhermore, if $u\in\mathcal{E}_m(\Om,\om)$ and $c={\rm constant}$ then $u+c\in\mathcal{E}_m(\Om,\om)$.	By Theorem \ref{thm 5.3}	we have 
		\begin{align*}
			&	( \omega + dd^c \max(u,v))^m\w\beta^{n-m} \\
			& \geq  \textit{1}_{\{u>v\}} ( \omega + dd^c \max(u,v))^m\w\beta^{n-m} +  \textit{1}_{\{u < v\}} ( \omega + dd^c \max(u,v))^m\w\beta^{n-m}\\
				& \geq  \textit{1}_{\{u>v\}} ( \omega + dd^c u)^m\w\beta^{n-m} +  \textit{1}_{\{u < v\}} ( \omega + dd^c v)^m\w\beta^{n-m}. 
		\end{align*}
		If $( \omega + dd^c v)^m\w\beta^{n-m} (\{u = v\}) = 0$, then the result follows. \\
Since $(\omega + dd^c v)^m\w\beta^{n-m}$ 
vanishes on $m$-polar sets,	the proof of \cite[Proposition 5.2]{HP17} shows that 
	$$ (\omega + dd^c v)^m\w\beta^{n-m}(\{ u = v + t \}) = 0, \; \;  \forall t \in \mathbb{R}\setminus I_{\mu}, $$
	where $I_{\mu}$ is at most countable. Take $\varepsilon_j \in \mathbb{R}\setminus I_{\mu}$,  $\varepsilon_j \searrow 0$.   We have 
	
	$$(\omega + dd^c v)^m\w\beta^{n-m}(\{ u = v + \varepsilon_j \}) = 0.$$
	This implies that
	
		\begin{align*} &( \omega + dd^c \max(u,v+ \varepsilon_j))^m\w\beta^{n-m} \\
			&\geq \textit{1}_{\{u> v+ \varepsilon_j \}} ( \omega + dd^c u)^m\w\beta^{n-m} + \textit{1}_{\{ u \leq v+ \varepsilon_j \}} ( \omega + dd^c v)^m\w\beta^{n-m}. 
			\end{align*}

\n	Let $\varepsilon_j \rightarrow 0$, by Theorem \ref{mono} together with using the Lebesgue Monotone Convergence Theorem  we get the desired conclusion.
	\end{proof}
Using the same arguments as in the proof of the above Corollary, we get the following Corollary.
	\begin{corollary}\label{cor 5.4}
	Let $u,v \in \mathcal{E}_m(\Omega,\omega)$, and let $\mu$ be a  positive Radon measure vanishing on $m$-polar sets. If 
	$$ 
	(\omega + dd^c u)^m\w\beta^{n-m} \geq \mu \; \emph{ and} \;  ( \omega + dd^c v)^m\w\beta^{n-m} \geq \mu,
	$$
	then
	$$  
	(\omega + dd^c \max(u,v))^m\w\beta^{n-m} \geq \mu. 
	$$
\end{corollary}

\section{ Envelope of $(\omega,m)$-subharmonic functions}
Motivated by \cite{ACCH09, HHHP14, HP17, Pmalay}, we will investigate the envelope of 
$(\omega,m)$--subharmonic functions and apply it to the study of Hessian equations in the next section.

Firstly, we need the following lemma.
\begin{lemma}\label{prop 2.1}
	Let $K \Subset \Omega$ be a relatively compact domain in $\Omega$. Fix $v \in \SH(\Omega,\omega)\cap L^\infty(\Omega)$ and $u \in \SH(K,\omega) \cap L^\infty(K)$. If $\limsup_{z\to\xi}u(z) \leq v(\xi)$ for all $\xi\in\partial K\cap\Omega$ then the function
	\[
	\varphi = 
	\begin{cases}
		\max(u,v) & \text{on } \; K \\
		v & \text{on} \; \Omega \setminus K, 
	\end{cases}
	\]
	belongs to $\SH(\Omega,\omega) \cap L^\infty(\Omega)$. \\
	Furthermore, if $\mu$ is a positive Radon measure on $\Omega$ such that 
	$$(\omega + dd^c v)^m\w\beta^{n-m} \geq \mu \; \text{on} \; \Omega \; \; \text{and} \; (\omega + dd^c u)^m\w\beta^{n-m} \geq \mu \; \text{on} \; K, $$
	then $(\omega + dd^c \varphi)^m\w\beta^{n-m} \geq \mu$ on $\Omega$.
\end{lemma}
\begin{proof}
	It is easy to see that $\varphi$ is upper semi-continuous on $\Omega$. 
	For $\varepsilon >0$, we consider
	\[
	\varphi_{\va} = 
	\begin{cases}
		\max(u,v+\va) & \text{on } \; K, \\
		v+\va & \text{on} \; \Omega \setminus K. 
	\end{cases}
	\]
	We will prove that $\varphi_{\va}\in \SH(\Omega,\omega)$ as below. \\
	\n (i) In this case $z\in K.$  There exists ball $B(z,r)\subset K$.Therefore, we infer that $\varphi_{\va}=\max(u,v+\va)$ is a $(\omega,m)$-sh function in $B(z,r).$\\
	\n (ii) Similarly, in this case $z\in \Omega\setminus \overline{K}.$ There exists a ball $B(z,r)\subset \Omega\setminus \overline{K}$. Hence, we have $\varphi_{\va}=v+\va$ is a $(\omega,m)$-sh function in $B(z,r).$\\
	\n (iii) In this case $z\in \partial K\cap\Omega.$ It follows from the hypothesis $\limsup_{z\to\xi}u(z) \leq v(\xi)$ $\forall$ $\xi\in\partial K\cap\Omega$ that there exists $r>0$ satisfying $\sup\limits_{B(z,r)}u\leq v(z)+\va.$ This implies that $\varphi_{\va}=v+\va$ in $B(z,r).$\\
	Combining (i), (ii) and (iii), we infer that $\varphi_{\va}$ is a $(\omega,m)$-sh function on $B(z,r)$ for all  $z\in\Omega.$ Therefore, $\varphi_{\va}$ is a $(\omega,m)$-sh function on $\Omega.$  Since $\varphi_{\va}\searrow \varphi$ as $\va\to 0,$ we obtain that $\varphi\in \SH(\Omega,\omega).$

	\n Next, we will prove the second assertion. Assume that $\mu$ is a positive Radon measure on $\Omega$ such that 
	$$(\omega + dd^c v)^m\w\beta^{n-m} \geq \mu \; \text{on} \; \Omega \; \; \text{and} \; (\omega + dd^c u)^m\w\beta^{n-m} \geq \mu \; \text{on} \; K. $$
	Note that, according to (ii) and (iii), we see that for all $z\in\Omega\setminus K,$ there exists a ball $B(z,r)$ such that $\varphi_{\va}=v+\va$ on $B(z,r).$ It implies that
	\begin{equation}\label{eqe7}(\omega + dd^c \varphi_{\va})^m\w\beta^{n-m}=\big(\omega + dd^c (v+\va)\big)^m\w\beta^{n-m}
	\end{equation}
	on $B(z,r).$ Hence, on $\Omega\setminus K$ we also have equality \eqref{eqe7}.\\

	According to Proposition  \ref{prop 5.2} and Remark \ref{rm5.3} (i), we get that $(\omega + dd^c v)^m\w\beta^{n-m} $ puts no mass on $m$-polar sets on $\Omega$. It follows from the hypothesis $(\omega + dd^c v)^m\w\beta^{n-m} \geq \mu \; \text{on} \; \Omega$ that $\mu$ puts no mass on $m$-polar sets on $\Omega.$ By Corollary \ref{cor 5.4} and the hypothesis, we infer that
	\begin{align*} &(\omega + dd^c \varphi_{\va})^m\w\beta^{n-m} = \ind_{K}(\omega + dd^c \varphi_{\va})^m\w\beta^{n-m}  + \ind_{\Omega\setminus K}(\omega + dd^c \varphi_{\va})^m\w\beta^{n-m} \\
		&  =  \ind_{K} \big(\omega + dd^c \max(u,v+\va)\big)^m\w\beta^{n-m} + \ind_{\Omega\setminus K}\big(\omega + dd^c (v+\va)\big)^m\w\beta^{n-m} \\
		&\geq  \ind_{K}\mu+ \ind_{\Omega\setminus K}\mu\\
		&\geq \mu.
	\end{align*}
	Letting $\va\searrow 0$, we get by Theorem \ref{mono} that 
	$$ (\omega + dd^c \varphi)^m\w\beta^{n-m} \geq \mu \; \text{on} \; \Omega.$$
	The proof is complete. 
\end{proof}

Now, we will prove the following result which is a generalization of \cite[Proposition 9.1]{BT82} and \cite[Lemma 4.1]{Pjmaa}.
\begin{proposition}\label{lem: constructing upper functions}
	Fix $\psi \in \SH(\Omega,\omega)\cap L^\infty(\Omega)$ with $\psi = 0$ in $\partial \Omega$. Assume $u \in \SH(\Omega,\omega) \cap L^\infty(\Omega)$ is such that $(\omega + dd^c u)^m\w\beta^{n-m} \geq (\omega + dd^c \psi)^m\w\beta^{n-m} $. Then, for every ball $B \Subset \Omega$, there is $u_{B} \in \SH(\Omega,\omega) \cap L^\infty(\Omega)$ such that
	\begin{center}
		$u_{B} \geq u$ on $\Omega,$ 
	\end{center}      
	\begin{center}
		$u_{B} = u$ on $\Omega \setminus B,$ 
	\end{center}  
	\begin{center}
		$(\omega + dd^c u_{B})^m\w\beta^{n-m} \geq (\omega + dd^c \psi)^m\w\beta^{n-m} ,$  
	\end{center}
	and
	\begin{center}
		$(\omega + dd^c u_{B})^m\w\beta^{n-m} = (\omega + dd^c \psi)^m\w\beta^{n-m} $  on $B$. 
	\end{center}
\end{proposition}
\begin{proof}
	The proof is divided into two steps.

	\textbf{Step 1.} Firstly, we  assume that there exists a function $h \in \mathcal{E}_m^0(B)$ such that
	$$ (\omega + dd^c \psi)^m\w\beta^{n-m} \leq (dd^c  h)^m\w\beta^{n-m}  \; \; \text{on} \; B, $$
	Let $u_j \in \mathcal{C}^0(\partial B)$ be a decreasing sequence that converges to $u$ on $\partial B$. By \cite[Lemma 9.3]{KN23c}, we can find a function $v_j \in \SH(B,\omega) \cap L^\infty(B)$ such that 
	\begin{equation}\label{eq8}v_j = u_j\,\,\text{ on}\,\, \partial B\end{equation} and
	\begin{equation}\label{eq7}(\omega + dd^c v_j)^m\w\beta^{n-m}= (\omega + dd^c \psi)^m\w\beta^{n-m}\,\,\text{ on}\,\, B.
	\end{equation}
	According to \cite[Corollary 3.11]{GN18}, we infer that the sequence $(v_j)$ is decreasing and $v_j \geq u$ on $B$. Assume that $ v_j\searrow v$ on $B$, then we have $v \in \SH(B,\omega) \cap L^\infty(B)$, $v\geq u$ on $B.$ It follows from equality \eqref{eq8} that $v = u$ on $\partial B.$ Moreover,  according to \eqref{eq7} and Theorem \ref{mono}, we infer that
	$$ (\omega + dd^c v)^m\w\beta^{n-m} = (\omega + dd^c \psi)^m\w\beta^{n-m} \,\,\text{ on}\,\,  B.$$
	
	Now, we	consider the function
	\[
	\varphi= 
	\begin{cases}
		v & \text{on } \; B, \\
		u & \text{on} \; \Omega \setminus B. 
	\end{cases}
	\]
	According to Proposition \ref{prop 2.1}, we know that $\varphi \in \SH(\Omega,\omega)\cap L^\infty$ and that  $$(\omega + dd^c \varphi)^m\w\beta^{n-m}\geq (\omega + dd^c \psi)^m\w\beta^{n-m}.$$ Now, we can take $u_B = \varphi$.

	\textbf{Step 2.} In general case, we fix $h \in \mathcal{E}_m^{0}(B)$ and set 
	$$ \mu_j = \ind_{B \cap \{\psi + \rho > j h\} \cup (\Omega \setminus B)} (\omega + dd^c \psi)^m\w\beta^{n-m}. $$
	Obviously, $\mu_j$ is a positive Radon measure on $\Omega$ and that $\mu_j \nearrow (\omega + dd^c \psi)^m\w\beta^{n-m}$ on $\Omega$ as $j \rightarrow +\infty$. 
	
	Observe that,  on $B$ we have
	\begin{equation}	\label{e6}
			\mu_j\leq \ind_{B \cap \{\psi + \rho > j h\} } (\omega + dd^c \psi)^m\w\beta^{n-m}\leq \ind_{B \cap \{\psi + \rho > j h\} } (dd^c (\psi+\rho))^m\w\beta^{n-m}. 
	\end{equation}
	Moreover, applying Corollary \ref{cor 5.5} in the case $\omega=0$, we get
	\begin{equation}
		\begin{split}\label{e7}
			&\ind_B\big( dd^c\max(\rho+\psi,jh)\big)^m\w\beta^{n-m}\\
			&\geq \ind_{B\cap\{\rho+\psi>jh\}}\big(dd^c(\rho+\psi)\big)^m\w\beta^{n-m}+\ind_{B\cap\{\rho+\psi\leq jh\}}\big(dd^c(jh)\big)^m\w\beta^{n-m}\\
			&\geq \ind_{B\cap\{\rho+\psi>jh\}}\big(dd^c(\rho+\psi)\big)^m\w\beta^{n-m}.
		\end{split}
	\end{equation}
	From inequality \eqref{e6} and inequality \eqref{e7}, we infer that
	$$ \mu_j \leq (dd^c \max(\rho + \psi, jh))^m\w\beta^{n-m} \; \; \textit{on} \; B. $$
	Since $\max(\rho + \psi, jh) \in \mathcal{E}_m^0(B)$, applying Step 1 with $(\omega + dd^c \psi)^m\w\beta^{n-m}$ is replaced by $\mu_j$, there exists a function 
	$u_{j,B} \in \SH(\Omega,\omega) \cap L^\infty$ which
	satisfies the following conditions:
	\begin{center}
		$u_{j,B} \geq u$ on $\Omega,$ 
	\end{center} 
	\begin{center}
		$u_{j,B} = u$ on $\Omega \setminus B,$ 
	\end{center}        
	
	\begin{center}
		$(\omega + dd^c u_{j,B})^m\w\beta^{n-m} \geq \mu_j$  on $\Omega$
	\end{center}
	and
	\begin{center}
		$(\omega + dd^c u_{j,B})^m\w\beta^{n-m} = \mu_j$  on $B$. 
	\end{center}
	Note that 
	\[
	u_{j,B}= 
	\begin{cases}
		v_j & \text{on } \; B, \\
		u & \text{on} \; \Omega \setminus B,
	\end{cases}
	\]
	where $v_j \in \SH(B,\omega)\cap L^{\infty}(B)$; $v_j\geq u$ on $B,$ $v_j=u$ on $\partial B$ and  $(\omega + dd^c v_j)^m\w\beta^{n-m}= \mu_j$ on $B.$
	Since $(\mu_j)$ is an increasing sequence of measures. By \cite[Corollary 3.11]{GN18} we infer that $(v_j)$ is a decreasing sequence of functions and so is $u_{j,B}$ as $j\to\infty.$ Assume that $u_{j,B}\searrow u_B$ as $j\to\infty.$ We see that $u_B$  has all the required properties. The proof is complete.
\end{proof}

\n Next, we study the envelope of $(\omega,m)$-sh functions.

\begin{theorem}\label{thm: omega-psh env}
	Let $v \in \mathcal{E}_m^0(\Omega) \cap \mathcal{C}^0(\bar\Omega)$, and let $\psi \in \SH(\Omega,\omega)\cap L^\infty$ with $\psi = 0$ on $\partial \Omega$. Then the function
	$$ u := \sup \{h \in \SH(\Omega,\omega) : (\omega+ dd^c h)^m\w\beta^{n-m} \geq (\omega+ dd^c \psi)^m\w\beta^{n-m}, \; \text{and} \; h \leq v - \rho\} $$
	belongs to $\SH(\Omega,\omega)\cap L^\infty(\Omega)$ and satisfies $$v+\psi \leq u \leq v - \rho \text{ on}\,\, \Omega,$$
	$$ (\omega + dd^c u)^m\w\beta^{n-m} \geq (\omega + dd^c \psi)^m\w\beta^{n-m} \; \; \text{on} \; \Omega ,$$ and
	$$ (\omega + dd^c u)^m\w\beta^{n-m} = (\omega + dd^c \psi)^m\w\beta^{n-m} \; \; \text{on} \; \{u < v-\rho\}.$$ 
\end{theorem}
\begin{proof}
	Observe that $v+\psi$ belongs to the class of functions that  construct $u.$ Hence, we get  $u^{*} \geq v+\psi$,  where $u^*$  is the upper semi-continuous regularization of $u.$. This implies that $u^{*} \in \SH(\Omega,\omega) \cap L^\infty(\Omega).$
	By Choquet's lemma, there is $u_j \in  \SH(\Omega,\omega) \cap L^\infty$ such that $u_j \leq v-\rho$, 
	\begin{equation}\label{e8}
		(\omega + dd^c u_j)^m\w\beta^{n-m} \geq  (\omega + dd^c \psi)^m\w\beta^{n-m}
	\end{equation} and $(\sup_j u_j)^* = u^*.$   By replacing $u_j$ with $\max(u_1,...,u_j)$,  according to Corollary \ref{cor 5.4}, we can assume that $(u_j)$ is non-decreasing. Therefore,  $u_j\nearrow u^*$   a.e.  on $\Omega$.  From $u_j \leq v-\rho$,  we infer that $u^* \leq v - \rho$ a.e on $\Omega.$ It implies that $u^*+\rho\leq v$ a.e on $\Omega.$  Since $u^*+\rho$ and $v$ are  $m$-sh functions, we obtain that $u^*+\rho\leq v$ everywhere on $\Omega.$ That is equivalent to $u^*\leq v-\rho$ on $\Omega.$ Moreover,  from Theorem \ref{mono} and inequality \eqref{e8}, we infer that 
$$(\omega + dd^c u^*)^m\w\beta^{n-m}  \geq (\omega + dd^c \psi)^m\w\beta^{n-m}\,\,\text{on}\,\,\Omega .$$
 By the definition of $u$, we obtain that $u^*\leq u.$ Hence $u = u^* \in \SH(\Omega,\omega) \cap L^\infty(\Omega),$ $$v+\psi\leq u\leq v-\rho$$ and 
	$$(\omega + dd^c u)^m\w\beta^{n-m}  \geq (\omega + dd^c \psi)^m\w\beta^{n-m} \,\,\text{on}\,\,\Omega  .$$
	
	It then remains to prove that $(\omega + dd^c u)^m\w\beta^{n-m} = (\omega + dd^c \psi)^m\w\beta^{n-m}$ on $\{u < v-\rho\}$.  Indeed, fix $z_0 \in \{u<v-\rho\}$. That means we have $u(z_0)<v(z_0)-\rho(z_0)$. It implies that $u(z_0)+\rho(z_0)-v(z_0)<0$. Since $v\in \mathcal{C}^0(\bar\Omega)$, we get $u+\rho-v$ is an upper semicontinuous function. Therefore, we can find a small ball $B$ centered at $z_0$ such that  $\sup\limits_{\bar{B}}(u+\rho-v)(z)<0$. Thus,  $\sup_{\bar B} (u+\rho) < \inf_{\bar B} v$. Therefore, we have $u(z)<v(z)-\rho(z)$ on $B.$ It suffices to prove $(\omega + dd^c u)^m\w\beta^{n-m} = (\omega + dd^c \psi)^m\w\beta^{n-m}$ on $B$. By Lemma \ref{lem: constructing upper functions}, there is $u_B \in \SH(\Omega,\omega) \cap L^\infty$ such that 
	\begin{equation}\label{e9}
		u_B \geq u \,\,\text{on}\,\, \Omega.
	\end{equation}
	\begin{equation}\label{e10}
		u_B = u\,\,\text{ on}\,\,\Omega \setminus B.
	\end{equation}
	
	\begin{equation}\label{e11}(\omega + dd^c u_B)^m\w\beta^{n-m} \geq (\omega + dd^c \psi)^m\w\beta^{n-m}\,\, \text{on}\,\, \Omega.
	\end{equation}
	\begin{equation}\label{e12}
		(\omega + dd^c u_B)^m\w\beta^{n-m} = (\omega + dd^c \psi)^m\w\beta^{n-m}\,\, \text{on}\,\, B.
	\end{equation}
	\n (i) Note that on $\Omega \setminus B,$ by equality \eqref{e10} we have  $u_B=u \leq v -\rho$.\\
	\n (ii) Moreover, since $u_B=u$ on $\partial B$ and $m$-sh functions $u_B+\rho, u+\rho$ achieves theirs maximum on the boundary, we have
	$$ \sup_{\bar B} (u_B+ \rho) = \sup_{\partial B} (u_B + \rho)= \sup_{\partial B} (u + \rho)=\sup_{\bar B} (u + \rho) < \inf_{\bar B} v.$$
	This implies that we  $u_B \leq v - \rho$ on $B.$\\
	Combining (i) and (ii) we infer that $u_B \leq v-\rho$ on $\Omega$. Coupling this with inequality \eqref{e11}, we see that $u_B$ belongs to the class of functions that construct $u.$ Thus $u_B\leq u.$.  It follows from inequality \eqref{e9} that $u_B=u.$ From equality \eqref{e12}, we get the desired. The proof is complete. 
\end{proof}

\n Now we will prove the main result in this section.
\begin{proposition}\label{main pro}
	Let $\psi \in \SH(\Omega,\omega) \cap L^\infty(\Omega)$ with $\psi = 0$ on $\partial \Omega$. Let $v \in \mathcal{E}_m(\Omega)$ be such that $(dd^c v)^m\w\beta^{n-m}$ is carried by a $m$-polar set. Then, the function 
	$$ \varphi = \sup \{h \in \mathcal{E}_m(\Omega,\omega) : (\omega + dd^c  h)^m\w\beta^{n-m} \geq (\omega + dd^c \psi)^m\w\beta^{n-m} \; \text{and} \; h \leq v -\rho\} $$
	belongs to $\mathcal{E}_m(\Omega,\omega)$, satisfies $\psi + v \leq \varphi \leq v - \rho$, and we have 
	\begin{equation}\label{e13} (\omega + dd^c \varphi)^m\w\beta^{n-m} = (\omega + dd^c \psi)^m\w\beta^{n-m} + (dd^cv)^m\w\beta^{n-m} \,\,\text{on}\,\,\Omega .\end{equation}
\end{proposition}

\begin{proof}
	Since $H_m(v)$ is carried by a $m$-polar set, we infer that 
	$$H_m(v)=\ind_{\{v=-\infty\}}H_m(v)+ \ind_{\{v>-\infty\}}H_m(v)= \ind_{\{v=-\infty\}}H_m(v).$$
	By \cite[Theorem 3.1]{Ch15}, there is $v_j \in \mathcal{E}_m^0(\Omega) \cap \mathcal{C}^0(\bar\Omega)$ such that $v_j \searrow v$. Consider
	$$ \varphi_j = \sup \{ h \in \mathcal{E}_m(\Omega,\omega) : (\omega + dd^c  h)^m\w\beta^{n-m}  \geq (\omega + dd^c \psi)^m\w\beta^{n-m}  \; \text{and} \; h \leq v_j -  \rho \}. $$ 
	By Theorem \ref{thm: omega-psh env}, we know that $\varphi_j \in \SH(\Omega,\omega) \cap L^\infty,$ 
	\begin{equation}\label{e14}
		v_j + \psi \leq \varphi_j \leq v_j - \rho,
	\end{equation}  
	\begin{equation} \label{e15} (\omega + dd^c \varphi_j)^m\w\beta^{n-m} \geq  (\omega + dd^c \psi)^m\w\beta^{n-m}
	\end{equation}
	and 	\begin{equation}\label{e16} \ind_{\{\varphi_j<v_j-\rho\}} (\omega + dd^c \varphi_j)^m\w\beta^{n-m} =  \ind_{\{\varphi_j<v_j-\rho\}}(\omega + dd^c \psi)^m\w\beta^{n-m}
	\end{equation}
	
	It follows from $v_j\searrow v$ that  $\varphi_j \searrow \varphi.$ Hence, we get $\varphi\in \SH(\Omega,\omega).$ Moreover,  from inequality \eqref{e14} we get   $v+\psi \leq \varphi \leq v-\rho$. This implies that $\varphi \in \mathcal{E}_m(\Omega,\omega).$\\
	\n Now, it remains to prove equality \eqref{e13}. Indeed, it follows from \eqref{e16} that
	\begin{equation}
		\begin{split}\label{e17}
			&(\omega + dd^c \varphi_j)^m\w\beta^{n-m}\\
			&=  \ind_{\{\varphi_j<v_j-\rho\}}(\omega + dd^c \varphi_j)^m\w\beta^{n-m}+\ind_{\{\varphi_j=v_j-\rho\}}(\omega + dd^c \varphi_j)^m\w\beta^{n-m}\\
			&=  \ind_{\{\varphi_j<v_j-\rho\}}(\omega + dd^c \psi)^m\w\beta^{n-m}+\ind_{\{\varphi_j+\rho=v_j\}}(\omega + dd^c \varphi_j)^m\w\beta^{n-m}\\
			&\leq  \ind_{\{\varphi_j<v_j-\rho\}}(\omega + dd^c \psi)^m\w\beta^{n-m}+ \ind_{\{\varphi_j+\rho=v_j\}}\big( dd^c (\rho+\varphi_j)\big)^m\w\beta^{n-m}\\
			&\leq (\omega + dd^c \psi)^m\w\beta^{n-m}+\ind_{\{\varphi_j+\rho=v_j\}}(dd^c v_j)^m\w\beta^{n-m}\\
			&\leq (\omega + dd^c \psi)^m\w\beta^{n-m}+ H_m(v_j).
		\end{split}
	\end{equation}
	where the second inequality follows from the fact that $\varphi_j+\rho\leq v_j$  and Corollary \ref{cor 5.5} in the special case when $\omega\equiv 0.$\\
	Combining inequality \eqref{e15} and inequality \eqref{e17} we infer that
	$$ (\omega + dd^c \psi)^m\w\beta^{n-m}\leq (\omega + dd^c \varphi_j)^m\w\beta^{n-m}\leq (\omega + dd^c \psi)^m\w\beta^{n-m}+ H_m(v_j)$$
	
	Taking $j\rightarrow +\infty$, by Theorem \ref{mono}, we infer that
	\begin{equation}\label{e18}
		(\omega + dd^c \psi)^m\w\beta^{n-m} \leq (\omega + dd^c \varphi)^m\w\beta^{n-m} \leq (\omega + dd^c \psi)^m\w\beta^{n-m} + H_m(v).    
	\end{equation}
	Moreover, since $\varphi\leq v-\rho,$ we infer that $\{\varphi>-\infty\}\subset \{v>-\infty\}.$ Therefore, we have
	$$ \ind_{\{\varphi>-\infty\}}H_m(v)\leq \ind_{\{v>-\infty\}}H_m(v)=0.$$ Hence, it follows from inequality \eqref{e18} that 
	\begin{equation}\label{e19}
		\ind_{\{\varphi>-\infty\}}(\omega + dd^c \psi)^m\w\beta^{n-m}= 	\ind_{\{\varphi>-\infty\}}(\omega + dd^c \varphi)^m\w\beta^{n-m}.
	\end{equation}
	On the other hand, we have $v + \psi \leq \varphi \leq v -\rho$. This implies that $\varphi +\rho\leq v.$ It follows from  Remark \ref{rm5.3} (ii) and \cite[Proposition 5.2]{HP17} that 
	\begin{align*} \ind_{\{\varphi=-\infty\}} (\omega + dd^c \varphi)^m\w\beta^{n-m} &= \ind_{\{\varphi=-\infty\}}\big(dd^c (\varphi+\rho)\big)^m\w\beta^{n-m}\\
		& =\ind_{\{\varphi+\rho=-\infty\}}\big(dd^c (\varphi+\rho)\big)^m\w\beta^{n-m}\\
		&\geq \ind_{\{v=-\infty\}} (dd^c v)^m\w\beta^{n-m} = H_m(v). 
	\end{align*}
	\n Moreover, from $v+\psi\leq\varphi,$  we infer that $v+\rho+\psi\leq \varphi+\rho.$ According to \cite[Proposition 5.2]{HP17}, we obtain that
	\begin{align*}
		&\ind_{\{\varphi=-\infty\}} (\omega + dd^c \varphi)^m\w\beta^{n-m}=	\ind_{\{\varphi=-\infty\}}\big(dd^c (\varphi+\rho)\big)^m\w\beta^{n-m}\\
		&=\ind_{\{\varphi+\rho=-\infty\}}\big(dd^c (\varphi+\rho)\big)^m\w\beta^{n-m}\\
		&\leq \ind_{\{v+\rho+\psi=-\infty\}}\big(dd^c (v+\rho+\psi)\big)^m\w\beta^{n-m}\\
		&=\ind_{\{v=-\infty\}}\big(dd^c (v+\rho+\psi)\big)^m\w\beta^{n-m}\\
		&\leq \ind_{\{v=-\infty\}}H_m(v)=H_m(v),
	\end{align*}
	where the last inequality follows from the fact that $\rho+\psi$ is a bounded $m$-subharmonic function and \cite[Lemma 5.6]{HP17}.\\
	Combining two above arguments, we have 
	$$\ind_{\{\varphi=-\infty\}} (\omega + dd^c \varphi)^m\w\beta^{n-m}=H_m(v). $$	
	From this and inequality \eqref{e19}, we obtain
	\begin{align*}
		(\omega + dd^c \varphi)^m\w\beta^{n-m} &= \ind_{\{\varphi>-\infty\}} (\omega + dd^c \varphi)^m\w\beta^{n-m} + \ind_{\{\varphi = -\infty\}} (\omega + dd^c \varphi)^m\w\beta^{n-m} \\
		&= \ind_{\{\varphi>-\infty\}} (\omega + dd^c \psi)^m\w\beta^{n-m} + H_m(v)\\
		&= (\omega + dd^c \psi)^m\w\beta^{n-m} + H_m(v),
	\end{align*}
	where the last equality follows from the fact that $ (\omega + dd^c \psi)^m\w\beta^{n-m}$ vanishes on all $m$-polar sets according to  Proposition \ref{prop 5.2} and Remark \ref{rm5.3}.
	The proof is complete. 
\end{proof}
\section{Degenerate Hessian equations}\label{sec: solving MAE}
In this Section, we will prove the main result of this paper.
\begin{theorem}\label{main thm}
	Assume $\mu \leq (dd^c u)^m\w\beta^{n-m} $ is dominated by the Hessian measure of a function $u \in \mathcal{K}_m(\Omega)\in\{\mathcal{E}_m(\Omega), \mathcal{N}_m(\Omega), \mathcal{F}_m(\Omega) \}$. Then, there is $\varphi \in \mathcal{K}_m(\Omega,\omega,\phi)$ such that 
	$$ (\omega + dd^c \varphi)^m\w\beta^{n-m}= \mu. $$
\end{theorem}
\begin{proof}
	According to \cite[Theorem 2.15]{Gasmi}, we can write 
	$$\mu = f(dd^c \psi)^m\w\beta^{n-m} +\nu, $$
	where $\psi \in \mathcal{E}_m^0(\Omega)$, $0 \leq f \in L^1\big((dd^c \psi)^m\w\beta^{n-m}\big)$ and $\nu$ is a positive Radon measure carried by a $m$-polar set. From hypothesis $\mu\leq H_m(u),$ we infer that $\nu\leq H_m(u).$ Applying \cite[main Theorem]{Gasmi}, there exists a function $v \in \mathcal{K}_m(\Omega,\phi)$ such that $\nu = (dd^c v)^m\w\beta^{n-m}.$ \\
	Put $\mu_j=\min(f,j) (dd^c \psi)^m\w\beta^{n-m} .$ We have $$\mu_j\leq jH_m(\psi)=H_m(\frac{\psi}{\sqrt[m]{j}}).$$ 
	By \cite[Lemma 9.3]{KN23c}, there is $\psi_j \in \SH(\Omega,\omega) \cap L^\infty$ such that $\psi_j = 0$ on $\partial \Omega$ and 
	\begin{equation}\label{e21}  (\omega + dd^c \psi_j)^m\w\beta^{n-m} =\mu_j= \min(f,j) (dd^c \psi)^m\w\beta^{n-m}. \end{equation}
	Moreover, \cite[Corollary 3.11]{GN18} shows that $\psi_j$ is decreasing. We assume that $  \psi_j\searrow \tilde{\psi}$. Then, we have $\tilde{\psi}\in \SH(\Omega,\omega).$\\
	 We will prove that $\tilde{\psi}\in \mathcal{K}_m(\Omega,\omega).$ Indeed, it follows from $\mu\leq H_m(u)$ and \cite[main Theorem]{Gasmi}, there exists a function $\xi\in\mathcal{K}_m(\Omega)$ such that $\mu=H_m(\xi)$ and $u+\phi\leq\xi\leq \phi.$ Choose $\rho\in C^2(\overline{\Omega})\cap \SH(\Omega)$ such that $\rho=0$ on $\partial \Omega$, $dd^c\rho\geq \omega$ and  $dd^c\rho\geq -\omega.$ Now, we have 
	$$(\omega+dd^c\psi_j)^m\w\beta^{n-m}=\mu_j\leq \mu=H_m(\xi)\leq \big(\omega+dd^c(\xi+\rho)\big)^m\w\beta^{n-m}.$$
	According to Proposition \ref{prop 5.2} and Remark \ref{rm5.3}, we have $(\omega+dd^c\psi_j)^m\w\beta^{n-m}$ vanishes on $m$-polar sets. By Corollary \ref{cor 5.4} we infer that
	$$(\omega+dd^c\psi_j)^m\w\beta^{n-m}\leq  \big(\omega+dd^c\max(\xi+\rho,\psi_j)\big)^m\w\beta^{n-m}.$$ 
	Moreover, we also have $\psi_j\geq \max(\xi+\rho,\psi_j)$ on $\partial \Omega.$ According to \cite[Corollary 3.11]{GN18}, we get $\psi_j\geq \xi+\rho$ on $\Omega$ for all $j\geq 1.$ Letting $j\to\infty,$ it implies that $\tilde{\psi}\geq \xi+\rho.$ Hence, $\rho+\tilde{\psi}\geq \xi+2\rho.$ This means $\tilde{\psi}\in \mathcal{K}_m(\Omega,\omega)$ as the desired.\\
	
	It follows from equality \eqref{e21} and Theorem \ref{mono} that
	$$  (\omega + dd^c \tilde{\psi})^m\w\beta^{n-m} = f (dd^c \psi)^m\w\beta^{n-m}. $$
	Consider 
	$$ \varphi_j := \sup \{h \in \mathcal{E}_m(\Omega,\omega) : (\omega + dd^ch)^m\w\beta^{n-m} \geq (\omega+ dd^c \psi_j)^m\w\beta^{n-m} \; \text{and} \; h \leq v - \rho \}.  $$
	Applying Proposition \ref{main pro}, with   $\psi_j$ playing the role of $\psi$, we infer that $\varphi_j \in \mathcal{E}_m(\Omega,\omega)$, 
	\begin{equation}\label{e22}\psi_j + v \leq \varphi_j \leq v- \rho
	\end{equation} and 
	\begin{equation}\label{eq MA meas varphi j}
		(\omega + dd^c \varphi_j)^m\w\beta^{n-m} = (\omega+ dd^c \psi_j)^m\w\beta^{n-m} + H_m(v).
	\end{equation}
	Moreover, since
	$$ (\omega+ dd^c \psi_j)^m\w\beta^{n-m} = \min(f,j) (dd^c\psi)^m\w\beta^{n-m}$$
	is increasing with respect to $j$, by the definition of $\varphi_j$, we infer that the sequence $(\varphi_j)$ is decreasing. Assume that $ \varphi_j\searrow \varphi $, then we have $\varphi \in \SH(\Omega,\omega).$ From inequality \eqref{e22}, we get  $\tilde{\psi} + v \leq \varphi \leq v- \rho$. Hence, we obtain that $$\tilde{\psi} + v+\rho\leq\varphi+\rho\leq v\leq\phi.$$ Note that, since $\tilde{\psi}\in\mathcal{K}_m(\Omega,\omega)$, we get $\tilde{\psi}+\rho\in \mathcal{K}_m(\Omega).$ Moreover, since $v\in \mathcal{K}_m(\Omega,\phi),$ there exists a function $\tilde{v}\in\mathcal{K}_m(\Omega)$ such that $\tilde{v}+\phi\leq v\leq \phi.$ Thus, we have $$\tilde{v}+\phi+\tilde{\psi}+\rho\leq v+ \tilde{\psi}+\rho\leq \varphi+\rho\leq v\leq \phi.$$  It implies that $\varphi \in \mathcal{K}_m(\Omega,\omega,\phi)$.  Letting $j \rightarrow +\infty$ in \eqref{eq MA meas varphi j},  by Theorem \ref{mono} we obtain that
	$$ (\omega + dd^c \varphi)^m\w\beta^{n-m} =  \mu, $$
	The proof is complete. 
\end{proof}
\section{A comparison principle}
Note that we do not know whether the solution of the Hessian equation in Theorem \ref{main thm} is unique. However, based on \cite[Theorem 5.5]{ACLR25b}, in this Section, we will prove a comparison principle. Firstly, based on the work of \cite{GZ07,ACLR25,ACLR25b}, we will define non-$m$-polar part of Hessian measure of $(\omega,m)$-subharmonic functions as below.

\n Fix $u \in \SH(\Omega,\omega)$. Fix $t$, $s\in \mathbb{R}_+$ and set $u_t = \max(u,-t)$. By Corollary \ref{cor 5.4}, we know that 
$$ \ind_{\{u_{t}>-s\}}(\omega + dd^c u_{t})^m\w\beta^{n-m} = \ind_{\{u_{t}>-s\}}(\omega + dd^c \max(u_{t},-s))^m\w\beta^{n-m}.  $$
For $t\geq s$, we have 
\begin{align*}
	\ind_{\{u>-t\}}(\omega + dd^c u_{t})^m\w\beta^{n-m} &\geq 
	\ind_{\{u>-s\}}(\omega + dd^c u_{t})^m\w\beta^{n-m} \\
	&=  \ind_{\{u_{t}>-s\}}(\omega + dd^c u_{t})^m\w\beta^{n-m} \\
	&=  \ind_{\{u_{t}>-s\}}\big(\omega + dd^c\max( u_{t},-s)\big)^m\w\beta^{n-m} \\
	&=  \ind_{\{u_{t}>-s\}}\big(\omega + dd^c\max( u,-s)\big)^m\w\beta^{n-m} \\
	&= \ind_{\{u>-s\}}(\omega + dd^c u_{s})^m\w\beta^{n-m}. 
\end{align*}
It implies that the sequence of general term
$$ \ind_{\{u> -t\}} (\omega + dd^c \max(u,-t))^\w\beta^{n-m}, \; \; t \geq 0, $$
is non-decreasing. Moreover, according to Proposition \ref{prop 5.2} and Remark \ref{rm5.3},  we see that each term puts no mass on $m$-polar set. Therefore, we define the non-$m$-polar Hessian measure of $u$ by the limit
$$ NP(\omega + dd^c u)^m\w\beta^{n-m} := \lim_{t\rightarrow +\infty} \ind_{\{u> -t\}} (\omega + dd^c \max(u,-t))^m\w\beta^{n-m}. $$
Obviously, the non-$m$-polar part of Hessian measure  $(\omega + dd^c u)^m\w\beta^{n-m}$ vanishes on all $m$-polar sets for every $(\omega,m)$-subharmonic function $u$.\\
In particular, if $u \in \mathcal{E}_m(\Omega,\omega)$ then  by Theorem \ref{thm 5.3} and Lebesgue monotone convergence theorem, we get
\begin{equation}
	\begin{split}\label{e3} NP(\omega + dd^c u)^m\w\beta^{n-m} &=\lim_{t\rightarrow +\infty} \ind_{\{u> -t\}} (\omega + dd^c \max(u,-t))^m\w\beta^{n-m}.\\
		&=\lim_{t\rightarrow +\infty} \ind_{\{u> -t\}} (\omega + dd^c u)^m\w\beta^{n-m}.\\
		&=\ind_{\{u> -\infty \}} (\omega + dd^c u)^m\w\beta^{n-m} 
	\end{split}
\end{equation} 
Furthermore, by Proposition \ref{pro3.10}, we obtain $\SH(\Omega,\omega)\cap L^{\infty}(\Omega)\subset \mathcal{E}_m(\Omega,\omega).$ Thus, if $u\in \SH(\Omega,\omega)\cap L^{\infty}(\Omega)$ then by Proposition \ref{prop 5.2} and Remark \ref{rm5.3}, we have
\begin{equation}\label{bounded}NP(\omega + dd^c u)^m\w\beta^{n-m}=(\omega + dd^c u)^m\w\beta^{n-m}
\end{equation}

The following result is a generalization of \cite[Theorem 2.2]{ACLR25}. 
\begin{theorem}\label{Dem ineq}
	Let $u$, $v \in \mathcal{E}_m(\Omega,\omega)$. We have 
	\begin{align*} NP(\omega + dd^c \max(u,v))^m\w\beta^{n-m} \geq  & \ind_{\{u \geq v\}} NP(\omega + dd^c u)^m\w\beta^{n-m}\\
		&  +  \ind_{\{u<v\}} NP(\omega + dd^c v)^m\w\beta^{n-m} . 
	\end{align*}
	In particular, if $u \leq v$ then 
	$$\ind_{\{u = v\}} NP(\omega + dd^c u)^m\w\beta^{n-m}  \leq \ind_{\{u = v\}} NP(\omega + dd^c v)^m\w\beta^{n-m} . $$
\end{theorem}
\begin{proof}
	For every $t \geq 0$, we put $u_t = \max(u,-t)$ and $v_t = \max(v,-t)$. Since $u_t$ and $v_t$ are bounded $(\omega,m)$-subharmonic functions on $\Omega$, according to equation \eqref{bounded}, we get $(\omega+dd^cu_t)^m\w\beta^{n-m}$ and  $(\omega+dd^cv_t)^m\w\beta^{n-m}$  put no mass on $m$-polar sets. By  Corollary \ref{cor 5.5} we have
	\begin{align*}  (\omega + dd^c \max(u_t,v_t))^m\w\beta^{n-m}& \geq  \ind_{\{u_t \geq v_t\}} (\omega + dd^c u_t)^m\w\beta^{n-m}\\ & +  \ind_{\{u_t<v_t\}} (\omega + dd^c v_t)^m\w\beta^{n-m}. \end{align*}
	Multiplying both sides of the above inequality by $\ind_{\{\min(u,v)>-t\}}$, we get
	\begin{equation}
		\begin{split}\label{eq3}
			& \ind_{\{\min(u,v)>-t\}}(\omega + dd^c \max(u_t,v_t))^m\w\beta^{n-m}\\&\geq \ind_{\{\min(u,v)>-t\}\cap \{u_t \geq v_t\}}(\omega + dd^c u_t)^m\w\beta^{n-m}\\
			&+\ind_{\{\min(u,v)>-t\}\cap \{u_t < v_t\}}(\omega + dd^c v_t)^m\w\beta^{n-m}.
		\end{split}
	\end{equation}
	By Theorem \ref{thm 5.3}, we infer that
	\begin{align*} 
		&\ind_{\{\min(u,v)>-t\}\cap \{u_t \geq v_t\}}(\omega + dd^c u_t)^m\w\beta^{n-m}\\
		&=\ind_{\{\min(u,v)>-t\}\cap \{u \geq v\}}(\omega + dd^c u)^m\w\beta^{n-m}\\
		&=\ind_{\{u\geq v>-t\}} (\omega + dd^c u)^m\w\beta^{n-m}. 
	\end{align*}  	
	Similarly, we have
	\begin{align*} 
		&\ind_{\{\min(u,v)>-t\}\cap \{u_t < v_t\}}(\omega + dd^c v_t)^m\w\beta^{n-m}\\
		&=\ind_{\{\min(u,v)>-t\}\cap \{u < v\}}(\omega + dd^c v)^m\w\beta^{n-m}\\
		&=\ind_{\{v> u>-t\}} (\omega + dd^c v)^m\w\beta^{n-m}. 
	\end{align*}  
	
	Moreover, we also have
	\begin{align*} 
		&\ind_{\{\min(u,v)>-t\}}(\omega + dd^c \max(u_t,v_t))^m\w\beta^{n-m}\\
		&=	\ind_{\{\min(u,v)>-t\}\cap \{\max(u,v)>-t\}}\Big(\omega + dd^c \max\big(\max(u,v),-t\big)\Big)^m\w\beta^{n-m}\\
		&=	\ind_{\{\min(u,v)>-t\}\cap \{\max(u,v)>-t\}}\Big(\omega + dd^c \max(u,v)\Big)^m\w\beta^{n-m}\\
		&\leq \ind_{ \{\max(u,v)>-t\}}\Big(\omega + dd^c \max(u,v)\Big)^m\w\beta^{n-m}
	\end{align*} 	
	
	Thus, it follows from inequality \eqref{eq3} that
	
	\begin{align*}
		\ind_{ \{\max(u,v)>-t\}} (\omega + dd^c \max(u,v))^m\w\beta^{n-m} &\geq   \ind_{\{u\geq v>-t\}} (\omega + dd^c u)^m\w\beta^{n-m} \\ &+  \ind_{\{-t <u<v\}} (\omega + dd^c v)^m\w\beta^{n-m}.
	\end{align*}
	Letting $t \rightarrow +\infty$, by equality \eqref{e3} we infer that
	\begin{align*} &NP(\omega + dd^c \max(u,v))^m\w\beta^{n-m}\\ & \geq   \ind_{\{u \geq v>-\infty\}} (\omega + dd^c u)^m\w\beta^{n-m}
		+  \ind_{\{-\infty <u<v\}} (\omega + dd^c v)^m\w\beta^{n-m} \\
		&=\ind_{\{u \geq v\}\cap \{v>-\infty\}}\ind_{\{u>-\infty\}} (\omega + dd^c u)^m\w\beta^{n-m}\\
		& +  \ind_{\{u<v\}\cap \{u>-\infty\}}\ind_{\{v>-\infty\}} (\omega + dd^c v)^m\w\beta^{n-m} \\
		&=\ind_{\{u \geq v\}\cap \{v>-\infty\}} NP (\omega + dd^c u)^m\w\beta^{n-m}
		+\ind_{\{u<v\}\cap \{u>-\infty\}} NP(\omega + dd^c v)^m\w\beta^{n-m}\\
		&=\ind_{\{u \geq v\}} NP (\omega + dd^c u)^m\w\beta^{n-m}
		+\ind_{\{u<v\}} NP(\omega + dd^c v)^m\w\beta^{n-m}.
	\end{align*} 
	The proof is complete.
	
\end{proof}
\n Now, according to \cite[Definition 5.1]{ACLR25b}, we define a relation $\preceq$ among  $(\omega,m)$-subharmonic functions. Let $u,v\in\SH(\Omega,\omega).$ We say that $u$ is more singular than $v$, and write $u\preceq v,$ if for any compact subset $K\Subset\Omega,$ there exists a constant $C_K$ such that $u\leq v+C_K$ on $K.$\\
The following result is an extension of Theorem 5.5 in \cite{ACLR25b}. Note that, in the case $\omega=0, m=n$ we obtain Theorem 3.6 in \cite{ACCH09} without assuming condition $\int_{\Omega}-w(dd^cu)<+\infty$ for some $w\in\mathcal{E}(\Omega)$ which is not identically 0.
\begin{theorem}\label{thm : uniqueness}
	Let $\phi \in \mathcal{E}_m(\Omega)\cap \MSH_m(\Omega)$, and let $u$, $v \in \mathcal{N}_m(\Omega,\omega,\phi)$ be such that $u \preceq v$. 
	If $(\omega + dd^c u)^m\w\beta^{n-m} \leq ( \omega + dd^c v)^m\w\beta^{n-m}$  then $u \geq v$. In particular,
	if $(\omega + dd^c u)^m\w\beta^{n-m} = ( \omega + dd^c v)^m\w\beta^{n-m}$  then $u = v$.
\end{theorem}
\begin{proof}
	Let $\tilde{u} \in \mathcal{N}_m(\Omega)$ be such that $\tilde{u} +\phi \leq u + \rho \leq \phi$.  Observe that $$u - v = u + \rho - (v+\rho) \geq \tilde{u}+\phi-v-\rho\geq \tilde{u}.$$
	Therefore, setting 
	$$ \varphi :=  \sup \{ h \in \SH(\Omega)\; ; \; h \leq \min(u - v,0) \}^*,  $$
	we have $\varphi \in \SH(\Omega)$ and $\tilde{u} \leq \varphi \leq 0$, hence $\varphi \in \mathcal{N}_m(\Omega)$.\\
	On one hand, it follows from the hypothesis $$(\omega + dd^c u)^m\w\beta^{n-m} \leq ( \omega + dd^c v)^m\w\beta^{n-m}$$ that
	\begin{equation}
		\begin{split}\label{e28}
			&\ind_{\{u=-\infty\}} (\omega + dd^c u)^m\w\beta^{n-m}\leq \ind_{\{u=-\infty\}}( \omega + dd^c v)^m\w\beta^{n-m}\\
			&=\ind_{\{u=-\infty\}}\ind_{\{v=-\infty\}}( \omega + dd^c v)^m\w\beta^{n-m}+\ind_{\{u=-\infty\}}\ind_{\{v>-\infty\}}( \omega + dd^c v)^m\w\beta^{n-m}\\
			&=\ind_{\{u=-\infty\}}\ind_{\{v=-\infty\}}( \omega + dd^c v)^m\w\beta^{n-m}\\
			&\leq \ind_{\{v=-\infty\}}( \omega + dd^c v)^m\w\beta^{n-m}.
		\end{split}
	\end{equation}
	On the other hand, fix a ball $B\Subset\Omega.$ It follows from the hypothesis $u \preceq v$ that $u\leq v +C_B$ on $B.$  According to Proposition \ref{prop 2.5}, we obtain that 
		\begin{equation*} \ind_{\{v=-\infty\}} (\omega + dd^c v)^m\w\beta^{n-m}\leq \ind_{\{u=-\infty\}} (\omega + dd^c u)^m\w\beta^{n-m} \,\,\text{on}\,\,B.  
	\end{equation*}
Since $B$ can be chosen arbitrarily and $\Omega$ is covered by finitely many balls, we obtain
	\begin{equation}\label{e29} \ind_{\{v=-\infty\}} (\omega + dd^c v)^m\w\beta^{n-m}\leq \ind_{\{u=-\infty\}} (\omega + dd^c u)^m\w\beta^{n-m} \,\,\text{on}\,\,\Omega.  
	\end{equation}
	Combining inequality \eqref{e28} and inequality \eqref{e29}, we infer that
	$$\ind_{\{v=-\infty\}} (\omega + dd^c v)^m\w\beta^{n-m}= \ind_{\{u=-\infty\}} (\omega + dd^c u)^m\w\beta^{n-m} $$
	Hence, by Remark \ref{rm5.3} (ii) we have
	$$ \ind_{\{u=-\infty\}} \big(dd^c(u+\rho)\big)^m\w\beta^{n-m} = \ind_{\{v=-\infty\}} \big(dd^c(v+\rho)\big)^m\w\beta^{n-m}.$$
Thus, it follows from \cite[Lemma 5.3]{ACLR25b} that $$\varphi=P(u-v,0)=P[(u+\rho)-(v+\rho),0]\in\mathcal{E}_m^a(\Omega).$$ Hence, $(dd^c\varphi)^m\w\beta^{n-m}$ does not charge $m$-polar sets.\\
 Moreover, by \cite[Lemma 3.1, (3)]{ACLR25b}, we infer that $$\mu_r(\varphi) =\ind_{\{\varphi>-\infty\}} (dd^c \varphi)^m\w\beta^{n-m}=(dd^c \varphi)^m\w\beta^{n-m}.$$
	
	Consider now $D = \{\varphi = \min(u - v,0)\}$. It follows from \cite[Theorem 4.2]{ACLR25b} that $(dd^c \varphi)^m\w\beta^{n-m}$ is carried by $D$. Now, we have
	\begin{equation}
		\begin{split}\label{e30}
			&(dd^c \varphi)^m\w\beta^{n-m}+ \ind_D NP(\omega + dd^c v)^m\w\beta^{n-m}\\
			&=(dd^c \varphi)^m\w\beta^{n-m}+ \ind_D \ind_{\{v>-\infty\}}(\omega + dd^c v)^m\w\beta^{n-m}\\
			&= (dd^c \varphi)^m\w\beta^{n-m} + \ind_{D \cap \{v>-\infty\}} (\omega + dd^c v)^m\w\beta^{n-m}\\
			&= \ind_{D \cap \{v>-\infty\}} \left( (dd^c \varphi)^m\w\beta^{n-m}+ (\omega + dd^c v)^m\w\beta^{n-m} \right) \\
			&=\ind_{D \cap \{v>-\infty\}\cap\{\varphi>-\infty\}} \left( (dd^c \varphi)^m\w\beta^{n-m}+ (\omega + dd^c v)^m\w\beta^{n-m} \right)\\
			&\leq\ind_{D \cap \{v>-\infty\}\cap\{\varphi>-\infty\}} \left(  (\omega + dd^c (v+\varphi) )^m\w\beta^{n-m} \right)\\
			&=\ind_{D \cap \{v>-\infty\}\cap\{\varphi>-\infty\}}\ind_{\{v+\varphi>-\infty\}} \left(  (\omega + dd^c (v+\varphi) )^m\w\beta^{n-m} \right)\\
			&\leq \ind_D\ind_{\{v+\varphi>-\infty\}}\left(  (\omega + dd^c (v+\varphi) )^m\w\beta^{n-m} \right)\\
			&=\ind_D NP\left(  (\omega + dd^c (v+\varphi) )^m\w\beta^{n-m} \right).
		\end{split}
	\end{equation}
	Since $\varphi + v \leq u$, by Theorem  \ref{Dem ineq}, we infer that
	$$
	\ind_{\{ \varphi + v = u\}}NP\left(  (\omega + dd^c (v+\varphi) )^m\w\beta^{n-m}\right)\leq  \ind_{\{ \varphi + v = u\}} NP(\omega + dd^c u)^m\w\beta^{n-m}.$$ Hence, we infer that
	\begin{align*}
		&\ind_{\{ \varphi + v = u\}\cap\{u< v\}}\ind_DNP\left(  (\omega + dd^c (v+\varphi) )^m\w\beta^{n-m}\right)\\
		&\leq \ind_{\{ \varphi + v = u\}\cap\{u< v\}}\ind_D NP(\omega + dd^c u)^m\w\beta^{n-m}.
	\end{align*}
	It implies that
	\begin{equation}\label{e31}
		\ind_{D}NP\left(  (\omega + dd^c (v+\varphi) )^m\w\beta^{n-m}\right)\leq  \ind_{D} NP(\omega + dd^c u)^m\w\beta^{n-m}.
	\end{equation}
	Moreover, it follows from the hypothesis $ (\omega + dd^c u)^m\w\beta^{n-m}\leq (\omega + dd^c v)^m\w\beta^{n-m}$ that
	\begin{align*} \ind_{\{u>-\infty\}} (\omega + dd^c u)^m\w\beta^{n-m}&\leq \ind_{\{u>-\infty\}}(\omega + dd^c v)^m\w\beta^{n-m}\\
		&\leq \ind_{\{v>-\infty\}}(\omega + dd^c v)^m\w\beta^{n-m}\\
		&=NP(\omega + dd^c v)^m\w\beta^{n-m}.
	\end{align*}
	That means we get
	\begin{equation}\label{e32}NP(\omega + dd^c u)^m\w\beta^{n-m}\leq NP(\omega + dd^c v)^m\w\beta^{n-m}.
	\end{equation}
	Combining inequality \eqref{e30}, inequality \eqref{e31} and inequality \eqref{e32}, we obtain that
	$(dd^c \varphi)^m\w\beta^{n-m} = 0$. It then follows from \cite[Proposition 5.4]{ACLR25b} that $\varphi = 0$, and thus $u \geq v$. 
	
	The second assertion follows directly from the first. Indeed, according to the first assertion we have $u \geq v.$ Thus, we have $v \preceq u.$ Changing the role of $u$ and $v$ in the first argument, we obtain that $v\leq u.$ Therefore, we obtain that $u = v$. The proof is complete.
\end{proof}

	\section*{Declarations}
	\subsection*{Ethical Approval}
	This declaration is not applicable.
	\subsection*{Competing interests}
	The authors have no conflicts of interest to declare that are relevant to the content of this article.
	\subsection*{Authors' contributions }
Nguyen Van Phu and Le Mau Hai  together studied  the manuscript.
	\subsection*{Funding }
	No funding was received for conducting this study.
	\subsection*{Availability of data and materials}
	This declaration is not applicable.


\begin{thebibliography}{000000}
		
	\bibitem[{\AA}CCH09]{ACCH09}
	P.~{\AA}hag, U.~Cegrell, R.~Czy{\.z}, and P. H. Hiep,
	\emph{{M}onge--{A}mp{è}re measures on pluripolar sets}, Journal de
	math{é}matiques pures et appliqu{é}es \textbf{92} (2009), no.~6,
	613--627.
	
	\bibitem[{\AA}CLR25]{ACLR25}
	P.~{\AA}hag, R.~Czy{\.z}, C.H. Lu, and A. Rashkovskii,  \emph{Geodesic connectivity and rooftop envelopes in the Cegrell classes}, Math. Ann. 391, 3333–3361 (2025). https://doi.org/10.1007/s00208-024-03003-7
	
	\bibitem[{\AA}CLR25b]{ACLR25b}
	P.~{\AA}hag, R.~Czy{\.z}, C.H. Lu, and A. Rashkovskii,  \emph{Kiselman minimum principle and rooftop envelopes in complex Hessian equations}, Math. Z. 308, article ID 70 (2024). 
	
	
	
	\bibitem[BT76]{BT76}
	E. Bedford and B. A. Taylor, \emph{The Dirichlet problem for a complex Monge-Ampère operator,}  Invent. Math. 37, 1–44 (1976). 
	
	\bibitem[BT82]{BT82}
	E. Bedford and B. A. Taylor, \emph{A new capacity for plurisubharmonic functions,} Acta. Math. 149, 1–40 (1982). 
	
	
	
	\bibitem[Bl05] {Bl05} Z. B{\l}ocki, {\it Weak solutions to the complex Hessian equation}, Ann. Inst. Fourier (Grenoble), {\bf 55}, 1735-1756 (2005).
	
	
	\bibitem[Ceg98]{Ceg98}
	U. Cegrell, \emph{Pluricomplex energy,} Acta Math, 180 (1998), 187-217.  DOI: 10.1007/BF02392899
	
	\bibitem[Ceg04]{Ceg04}
	U. Cegrell, \emph{The general definition of the complex Monge-Ampère operator,} Ann. Inst. Fourier (Grenoble) 54 (2004), 159-179. (www.numdam.org/articles/10.5802/aif.2014/)
	
	
	
	
	
	\bibitem[Ch12]{Ch12} L. H. Chinh, {\it On Cegrell's classes of $m-$subharmonic functions}, arXiv 1301.6502v1.
	
	\bibitem[Ch15]{Ch15} L. H. Chinh, {\it A variational approach to complex Hessian equation in $\mathbb{C}^{n}$}, J. Math. Anal. Appl., {\bf 431}, no. 1, 228-259 (2015).
	
	\bibitem [Cu13]{Cu13} N. N. Cuong, {\it Subsolution theorem for the complex Hessian equation}, Univ. Iagel. Acta Math.,  {\bf 50}, 69-88  (2013).
	
		\bibitem[Dem91]{Dem91}
	J. P. Demailly, \emph{Potential Theory in Several Complex Variables,}
	(http://www-fourier.ujf-grenoble.fr/$^{\sim}$demailly/books.html), (2016).
	
	\bibitem[Ga21]{Gasmi} A.E Gasmi, {\it The Dirichlet problem for the complex Hessian operator in the class $\mathcal N_m (\Omega, f)$}, Math. Scand., {\bf 121} (2021), 287-316. 
	
	
	\bibitem [GN18]{GN18} D. Gu, N.-C. Nguyen, {\it The Dirichlet problem for a complex Hessian equation on compact Hermitian manifolds with boundary}, Ann. Sc. Norm. Super. Pisa Cl. Sci., {\bf 18} (2018), no. 4, 1189-1248.
	
	
	
	
	\bibitem[GZ07]{GZ07} V. Guedj and A. Zeriahi, {\it The weighted Monge - Ampere energy of quasi-plurisubharmonic functions}, J. Funct. Anal. 250 (2) (2007), 442-482.
	
	
	
	
	\bibitem[HHHP14] {HHHP14} L. M. Hai, P. H. Hiep, N. X. Hong and N. V. Phu, {\it The Monge-Amp\`ere type equation in the weighted pluricomplex energy class}, Int. J. Math., {\bf 25}, no. 05, article ID: 1450042 (2014).
	
	
	\bibitem[HP17]{HP17} V. V. Hung and N. V. Phu, {\it Hessian measures on m- polar sets and applications to the complex Hessian equations}, Complex Var. Elliptic Equ., {\bf 62} (2017), no. 8, 1135-1164.
	
		\bibitem [Kl91]{Kl91} M. Klimek, Pluripotential Theory, Oxford Science Publications, The Clarendon Press Oxford University Press, New York, 1991.
		
	\bibitem[Ko95]{Ko95}
	S. Ko{\l}odziej, \emph{The range of the complex Monge-Ampère operator II},
	Indiana Univ. Math. J. 44 (1995), 765-782.
	
	
	
	\bibitem[KN15]{KN15Phong}
	S. Ko{\l}odziej and N. C. Nguyen, \emph{Weak solutions to the complex Monge-Ampère equation on Hermitian manifolds,} Analysis, complex geometry, and mathematical physics.  Contemporary Mathematics, vol. 644 (American Mathematical Society, Providence, RI, 2015) 141-158. 
	
	
	
	\bibitem[KN23a]{KN23a}
	S. Ko{\l}odziej and N. C. Nguyen, \emph{The Dirichlet problem for the Monge-Ampère equation on Hermitian manifolds with boundary,} Calc. Var. Partial Differential Equations
	62 (2023), no. 1, Paper No. 1.
	
	\bibitem[KN23b]{KN23b}
	S. Ko\l odziej and N. C. Nguyen, \emph{Weak Solutions to Monge–Ampère Type Equations on Compact Hermitian Manifold with Boundary,} The Journal of Geometric Analysis 33.1 (2023): 15.
	
	\bibitem[KN23c] {KN23c} S. Ko{\l}odziej and N.-C. Nguyen, {\it Complex Hessian measures with respect to a background Hermitian form}, https://arxiv.org/abs/2308.10405, to appear in Analysis and PDE.
	
	\bibitem [PD2023] {PD2023}  N. V. Phu and N. Q. Dieu, {\it Complex $m-$Hessian type equations in $\mathcal{E}_{m,\chi}(\Omega)$}, Publicationes Mathematicae Debrecen, {\bf 106}, no. 1-2 (13), 241-263 (2025).
	
	\bibitem[PDtaiwan] {PDtaiwan}  N. V. Phu and N. Q. Dieu, {\it Complex $m-$Hessian type equations in weighted energy classes of $m$-subharmonic functions with give boundary value},  Taiwanese Journal of Mathematics, {\bf 29}, no. 2, 413-424 (2025).
	
	
	\bibitem[Pjmaa]{Pjmaa}  N. V. Phu, {\it Approximation of $m$-subharmonic function with given boundary values}, J. Math. Anal. Appl., {\bf 534}, no. 2, 128097 (2024).
	
	\bibitem[Pmalay]{Pmalay}  N. V. Phu, {\it Approximation of $m$-subharmonic functions in weighted energy classes with given boundary values}, Bull. Malays. Math. Soc., {\bf 48} (2025), article number: 89.
	
	
	
	\bibitem[SA12] {SA12} A. S. Sadullaev and B. I. Abdullaev, {\it Potential theory in the class of $m-$subharmonic functions}, Proc. Steklov Inst. Math., {\bf 279}, 155-180 (2012).
	
	
	\bibitem[Sal25]{Sal25}
	M. Salouf, \emph{Degenerate complex Monge-Ampère equations with non-K\"ahler forms in bounded domains}, Indiana Univ. Math. J. 74 No. 1 (2025), 131–156.
	
	\bibitem[T19]{T19} N. V. Thien, {\it Maximal $m$-subharmonic functions and the Cegrell class $\mathcal{N}_{m}$ }, Indagationes Mathematicae, {\bf 30} (2019), Issue 4, 717-739.
	
	
	
	
	
	
	
	
		
		
		
	
		
		
	
		
		
		
		
		
	
		
		
		
		
	
		
		
		
	
	
		
		
		
		
		
	
		
	
	
		
		
			
		
		
		
	
		
		
	

		
	
		
		
		
		
	\end{thebibliography}
\end{document}